\documentclass{article}
\usepackage{amsmath}
\usepackage{amsthm}
\usepackage{amsfonts}
\usepackage{amssymb}
\usepackage{mathtools}
\usepackage[numbers]{natbib}
\usepackage{appendix}
\usepackage{xcolor}
\usepackage{tikz}
\usepackage{caption}
\usepackage{subcaption}
\usepackage[ruled,vlined,linesnumbered]{algorithm2e}

\newtheorem{theorem}{Theorem}
\newtheorem{definition}{Definition}

\newtheorem{proposition}{Proposition}

\newcommand{\N}{\mathbb N}
\newcommand{\E}{\mathbb E}
\newcommand{\R}{\mathbb R}
\newcommand{\M}{\mathcal M}
\newcommand{\F}{\mathcal F}
\newcommand{\norm}[1]{\left\lVert#1\right\rVert}
\newcommand{\X}{\mathcal{X}}
\newcommand{\Y}{\mathcal{Y}}

\title{Dual Stochastic Natural Gradient Descent
\\ 
\small{and convergence of interior half-space gradient approximations}}

\author{Borja S\'ANCHEZ-L\'OPEZ \and Jesús CERQUIDES\\ \\
\small{IIIA-CSIC, Campus UAB, Bellaterra, Spain}}
\date{}
\begin{document}
\maketitle
\begin{abstract}
The multinomial logistic regression (MLR) model is widely used in statistics and machine learning. Stochastic gradient descent (SGD) is the most common approach for determining the parameters of a MLR model in big data scenarios. However, SGD has slow sub-linear rates of convergence \cite{nemirovski_robust_2009}. A way to improve these rates of convergence is to use manifold optimization \cite{hu_brief_2020}. Along this line, stochastic natural gradient descent (SNGD), proposed by Amari \cite{amari_natural_1998}, was proven to be Fisher efficient when it converged. However, SNGD is not guaranteed to converge and it is computationally too expensive for MLR models with a large number of parameters. 

Here, we propose a stochastic optimization method for MLR based on manifold optimization concepts which (i) has per-iteration computational complexity is linear in the number of parameters and (ii) can be proven to converge. 

To achieve (i) we establish that the family of joint distributions for MLR is a dually flat manifold and we use that to speed up calculations. Sánchez-López and Cerquides \cite{borja_jesus_2019} have recently introduced convergent stochastic natural gradient descent (CSNGD), a variant of SNGD whose convergence is guaranteed. To obtain (ii) our algorithm uses the fundamental idea from CSNGD, thus relying on an independent sequence to build a bounded approximation of the natural gradient. We call the resulting algorithm dual stochastic natural gradient descent (DNSGD). By generalizing a result from Sunehag et al. \cite{sunehag_variable_2009}, we prove that DSNGD converges. Furthermore, we prove that the computational complexity of DSNGD iterations are linear on the number of variables of the model.
\end{abstract}

\section{Introduction}\label{sec:intro}
Multinomial Logistic Regression (MLR) is a  widely used tool for classification. Some relevant examples solving  real-world tasks are \cite{li_spectral_2012} for the image classification branch, \cite{covington_deep_nn_2016} for video recommendation tasks, or numerous examples in health and life sciences for analyzing nominal qualitative response variables, to name some \cite{daniels1997hierarchical, bull2007confidence, biesheuvel2008polytomous, leppink2020multicategory}.

The justification of MLR goes beyond practical. In statistical decision theory, it is well known that MLR for the choice probability can be derived assuming that (i) the random utilities are independent and identical distributed (i.i.d.) across alternatives and that (ii) their common distribution is a Gumbel function \cite{ben-akiva_discrete_1985}. Recent results \cite{tadei_recent_2018} show that the Gumbel distribution for the choice variables is not necessary and that any distribution which is asymptotically exponential in its tail is sufficient to obtain the MLR model.

Cross-entropy, also known as log-loss, is the loss function most used in MLR and it is also convenient from a statistical decision theory standout. Once the form of the loss function is elicited \cite{NIPS2008_077e29b1,reidCompositeBinaryLosses2010} and the inverse link function is understood as a mapping from scores to class probabilities, the log-loss is proved to be a proper composite loss together with logistic regression models \cite{reidCompositeBinaryLosses2010,vernetCompositeMulticlassLosses2011,nockSupervisedLearningNo2020}. This provides theoretical support for the usage of the log-loss from the standpoint of statistical decision theory. 


Classification algorithms predict the value of a discrete variable (class) given some other variables (features). We use $\Y$ for the class variable and $\X\in\Omega$ for the features. We assume a finite set of classes $\Y\in\{1,...,s\}$. We are interested in computing the unknown conditional probability distributions $P(\Y\mid \X)$. This is accomplished by optimizing the expected risk function \cite{vapnik_principles_1991}. Stochastic gradient descent (SGD), even though it was introduced in the mid-twentieth century, is the most common approach for the task because it is fast and simple. The strategy of SGD is very intuitive:  Assume $P$ is an unknown probability distribution on $\X\times\Y$ , $\mathbb{L}$ is a differentiable function from $\R^k$ to $\R$, and $l$ is a differentiable loss function such that $\mathbb{L}(\eta) = \E_{z\sim P}\left[l(\eta, z)\right]$. If $\eta_0\in\R^k$ is an initial estimation, SGD algorithm is defined by the update equation
\begin{equation}\label{eq:sgd}
    \eta_{t+1}\leftarrow \eta_t -\gamma_t\nabla l(\eta_t, s_t)
\end{equation}
where $\nabla l(\eta_t, s_t)=\frac{1}{|s_t|} \sum_{z\in s_t} \nabla l(\eta_t, z)$ is an approximation to the gradient of  $\mathbb{L}(\eta_t)$ according to a sample set $s_t\sim P$ and $\gamma_t$ is a positive number encoding the learning rate. After certain regularities on the function and the learning rate, SGD converges to the minimum \cite{bottou-98x}. However, usually convergence speed suddenly drops after a moderate solution quality has been reached, making the minimum unreachable from a practical point of view. Furthermore, SGD highly depends on the learning rate parameter, which can be difficult to tune, and it is vulnerable to the plateau phenomenon and to ill-conditioning. To face this issue, many SGD variants have arisen that basically modify the direction $\nabla l(\eta_t, s_t)$ using a positive semi-definite matrix $M_t$. Specifically, such generalization of SGD can be described by equation
\begin{equation}\label{eq:matrix_sgd}
    \eta_{t+1} \leftarrow \eta_t - \gamma_t M_t \nabla l(\eta_t, s_t)
\end{equation}
where $M_{t}$ is commonly a preconditioning matrix capturing the local curvature or related information such as the Hessian matrix in Newton´s method \cite{dennis_numerical_1996} or the inverse of the Fisher Information Matrix in Stochastic Natural Gradient Descent (SNGD) \cite{amari_natural_1998}. Due to the increment of the computational complexity (equation \ref{eq:matrix_sgd} defines a second order method) a trade-off between quality of curvature information and computational cost is assumed. Some widely used examples are preconditioned SGD, diagonal approximations of the Hessian \cite{becker_improving_1989}; Adagrad \cite{duchi_adaptive_2011}, Adadelta \cite{zeiler_adadelta:_2012}, RMSProp \cite{tieleman_lecture_2012} or Adam \cite{kingma_adam:_2015} that use the diagonal of the covariance matrix of the gradients.

Among existing preconditioning matrix algorithms, we focus our attention to SNGD and its variants. This kind of algorithm runs over a smooth manifold $\M$ of dimension $n$ \cite{carmo_riemannian_2013} equipped with a metric $g$ defined at every $p\in \M$. Metric $g_p$ is the positive-definite tensor that expresses the local metric information at $p$. The pair $(\M,g)$ is a Riemannian manifold of dimension $n$ \cite{carmo_riemannian_2013}. It is possible to choose a system of coordinates $\eta\in U\subset\R^n$ -- or parametrization -- to refer to points in $\M$.
In this case the metric information of $g$ at $\eta$ is given by an $n$ dimensional square matrix $G_\eta$, symmetric and positive-definite, in the base derived by the parametrization. 

Assume $\M$ is not standard $\R^n$, for instance when $\M$ is a sphere of dimension $n$ or when it is a Statistical manifold \cite{murray1993differential}, that is, manifolds whose points refer to probability distributions of the same family, where usually the Fisher information metric (FIM) is choosed. In such cases, it is interesting to work in a Riemannian manifold for two reasons. First, because it provides correct notions of angles  and local lengths which yields better updates. This is important when the gradient is the key tool for an optimization algorithm, since the gradient possesses these local magnitudes. And second, it allows to correctly define a direction in the space, in opposition to what happens with SGD, where the descent direction is parametrization dependent, in the sense that the gradient depends on the selected parametrization $\eta$.

Roughly speaking, the gradient of a function $f$ at $p\in\M$ is the steepest direction of $f$ at $p$. In a Riemannian manifold this is called the natural gradient by \cite{amari_natural_1998}, noted as $\widetilde{\nabla} f(p)$, and it is well defined since it takes into account the metric. If a parametrization $\eta$ is fixed, Amari \cite{amari_natural_1998} proved that
\begin{equation}\label{eq:natural_gradient}
    \widetilde{\nabla} f(\eta)=(G_\eta)^{-1}\nabla f(\eta)
\end{equation}
where $\widetilde{\nabla} f(\eta)$ is the natural gradient at $\eta$ in the base derived by the parametrization.

SNGD follows the natural gradient instead of the gradient. In particular, it sets the preconditioning matrix $M_t = (G_{\eta_t})^{-1}$ in update equation \ref{eq:matrix_sgd}, to follow the natural gradient according to equation \ref{eq:natural_gradient}
\begin{equation}
    \label{eq:sngd}
    \begin{split}
     \eta_{t+1} &\leftarrow \eta_{t} - \gamma_t \widetilde{\nabla} l(\eta_{t},z_t)\\
     &\leftarrow \eta_{t} - \gamma_t (G_{\eta_t})^{-1}{\nabla} l(\eta_{t},z_t)
     \end{split}
\end{equation} 
This algorithm, or an approximation, usually speeds up  learning in many problems, avoids the plateau effect and it defines parametrization independent directions. However, it faces two main problems: 

\begin{enumerate}
\item its computational complexity is high due to the need of either inverting a matrix or solving a linear system, and

\item it does not converge in some scenarios where SGD does \cite{thomas_2014} or it needs stronger assumptions such as compactness to stabilize \cite{bonnabel_stochastic_2013}.
\end{enumerate}
Issue i) warns about the higher computational complexity order, which really is a problem for nowadays large-scale high-dimensional problems, and issue ii) refers to the convergence property.

The objective of this work is to propose a natural gradient optimization method \cite{hu_brief_2020} for MLR, the Dual Stochastic Natural Gradient Descent (DSNGD), whose convergence is proven and whose computational complexity order equals that of SGD. Therefore this article reveals a strategy to approximate SNGD without suffering issues i) and ii). 

To deal with issue i) we establish that the family of joint distributions for MLR is a dually flat manifold and we use that to speed up calculations. To overcome issue ii) our algorithm uses the fundamental idea from CSNGD \cite{borja_jesus_2019}, relying on an independent sequence to build a bounded approximation of the natural gradient. 
 
Section \ref{sec:related_work} introduces some related work with essential concepts needed to solve issues i) and ii) for our natural gradient based algorithm. Section \ref{sec:dsngd} defines DSNGD. Sections \ref{sec:dsngd_complexity} and \ref{sec:dsngd_convergence} face issues i) and ii) respectively for the discrete case, that is, when $\X=\{1,...m\}$, and they prove discrete DSNGD is convergent and as fast as SGD, in terms of complexity order.

\section{Related work}\label{sec:related_work}

\subsection{Dually flat manifolds}\label{dfm}

Issue i) is adressed in this paper by restricting to dually flat manifolds (DFM) \cite{amari_information_2016,nielsen_elementary_2018}. 
The computational cost of natural gradient can be significantly reduced if the ambient space is a dually flat manifold, as one can see for instance for mirror descent \cite{nemirovski1983}.

DFM are built after two dual connections -- conjugate connections -- that are flat, that is, where Riemann-Christoffel curvature vanishes. As it is proved in \cite{amari_information_2016}, in such manifolds, there exist two dual parametrizations $\eta$ and $\eta^*$, related by the Legendre transform of a convex function $F(\eta)$, such that
\begin{equation}\label{eq:F_manifold}
    \begin{split}
        \eta^* = \nabla F(\eta)\\
        \nabla^2 F(\eta) = G_{\eta}
    \end{split}
\end{equation}
considering that $\eta$ and $\eta^*$ refer to the same point. This leads to a key property of DFM, starting by applying the chain rule to equation \ref{eq:natural_gradient}
\begin{equation}\label{eq:natural_gradient_is_gradient}
\begin{split}
    \widetilde{\nabla} f(\eta) &=(G_{\eta})^{-1}\nabla f(\eta)\\
    &=(G_{\eta})^{-1}\nabla \eta^*(\eta)\nabla f(\eta^*)\\
    &=(G_{\eta})^{-1}\nabla \nabla F(\eta)\nabla f(\eta^*)\\
    &=(G_{\eta})^{-1}G_\eta\nabla f(\eta^*)\\
    &=\nabla f(\eta^*)\\
\end{split}
\end{equation}
for any differentiable function $f$ defined in $\M$. 
Equation \ref{eq:natural_gradient_is_gradient} is proved for linear exponential families, a well known DFM, in
\cite{masegosa_stochastic_2014} and \cite{raskutti_information_2015}. Therefore, equation \ref{eq:natural_gradient_is_gradient} states that the natural gradient equals the gradient in its dual parametrization. Here one deduces a strategy to compute the natural gradient, or an approximation, without paying the costs of matrix inversion or linear system solving. Summing up, the main idea that solves issue ii) is equation \ref{eq:natural_gradient_is_gradient}. For example, in the case of SNGD in a DFM, one equivalently writes SNGD as
\begin{equation}
    \label{eq:fast_sngd}
     \eta_{t+1} \leftarrow \eta_{t} - \gamma_t \nabla l(\eta^*_t, z_t)
\end{equation}
Equations \ref{eq:sngd} and \ref{eq:fast_sngd} define SNGD, but the latter avoids matrix inversions and linear system solving. 
\subsection{Mirror descent}
DSNGD described in this article is by no means the only algorithm taking profit from dual space properties. Mirror descent algorithm \cite{nemirovski1983} makes use of dual parametrizations. As proved in \cite{raskutti_information_2015}, mirror descent in a dually flat manifold is nothing else than SNGD run in the dual space. 

According to \cite{BECK2003167}, mirror descent follows below update rule.

\begin{equation}
\begin{split}
    \eta_t \leftarrow& \nabla F^*(\eta^*_t)\\ 
    \eta^*_{t+1} \leftarrow& \nabla F(\eta_t) -\gamma_t \nabla l(\eta_t, s_t)\\ 
    \end{split}
\end{equation}
where $F$ is a convex function and $F^*$ is the Legendre transform of $F$. 
Even though both DSNGD and mirror descent rely on duality, there are clear differences between them: 
(i) mirror descent is normally defined for off-line learning, (ii) DSNGD has its convergence guaranteed, and (iii) mirror descent keeps a sequence of points in the manifold expressed both in primal and dual parametrizations while DSNGD has two sequences moving in the primal and dual space which are not necessarily connected. Difference (iii) is specially relevant since it requires mirror descent to rely on the computation of duals, while DSNGD can be run in spaces where we do not even know how to efficiently compute the dual coordinates of a point when given its primal coordinates.  

\subsection{Multinomial logistic regression}

As described above, our strategy to reduce the per-iteration computational complexity of DSNGD relies on the fact that the family of joint distributions for MLR is a dually flat manifold. 



The main assumption of MLR \cite{banerjee_analysis_2007} is that the log-odds ratio of the class posteriors $P(\Y\mid\X)$ is an affine function of the features $\X$. 
Banerjee \cite{banerjee_analysis_2007} proved (Theorem 2) that a class of distributions fulfills that assumption if and only if for each value of $\Y$, the class of conditional distributions $P(\X\mid \Y)$ belongs to the same linear exponential family (LEF)\footnote{For a definition of linear exponential family see \cite{wani_1968}}. 
In section \ref{sec:joint} we use these results to prove that the class of joint distributions $P(\Y,\X)$ is also a LEF. It is well known that a LEF is a DFM \cite{amari_information_2016}. Usually, finding the minimum expected risk MLR parameters is formulated as an optimization problem in $\mathbb{R}^k$ which is solved by means of SGD. Instead, we propose to formulate the problem as a manifold optimization problem \cite{hu_brief_2020}, over the manifold of probability distributions $P(\Y,\X)$ fulfilling the main assumption of MLR. Since we will prove that this manifold is dually flat, this formulation of the problem will allow us to capture the curvature information of the manifold efficiently.

\subsection{Convergent Stochastic Natural Gradient Descent}
In \cite{borja_jesus_2019}, an convergent variant for SNGD, namely CSNGD, is presented. CSNGD becomes stable in every toy scenario presented, unlike SNGD which fails in those same situations. Moreover, from a practical point of view it inherits the convergence speed of SNGD. CSNGD is defined with the update rule
\begin{equation}\label{eq:csngd}
    \eta_{t+1} \leftarrow \eta_t-\gamma_t (G_{\zeta_t})^{-1} {\nabla} l(\eta_{t},z_t)
\end{equation} 
where $\{\zeta_t\}_{t\in\N}$ can be any convergent sequence in $\R^k$. Both SNGD and CSNGD work by progressively building a sequence $\{\eta_t\}_{t\in\N}$. However, CSNGD additionally maintains an independent  sequence ${\zeta}_t$, which is only required  to be convergent. This difference allows CSNGD to converge to the unique minimum of a convex function after some reasonable conditions on the learning rate parameter. Precisely, SNGD does not converge due to the inverse matrix $(G_{\eta_t})^{-1}$ in equation \ref{eq:sngd}, since eigenvalues of that matrix are unbounded. CSNGD forces convergence and eigenvalue confinement of sequence $\{(G_{\zeta_t})^{-1}\}_{t\in\N}$ because of the convergence of sequence $\{\zeta_t\}_{t\in\N}$ and continuity property, stabilizing the algorithm. This idea is used in section~\ref{sec:dsngd} to define DSNGD and it allows us to prove its convergence in section~\ref{sec:dsngd_convergence}.

\subsection{Convergence of interior half-space gradient approximations}

In \cite{sunehag_variable_2009} Sunehag et al. provide a variable metric stochastic approximation theory. One of the key results in that paper is Theorem 3.2 which proves convergence given that we take a scaling matrix $B_t$ at step $t$ of the algorithm, provided that the spectrum of their (possibly non-convergent) scaling matrices is uniformly bounded from above by a finite constant and from below by a strictly positive constant. Moreover it assumes the step direction at iteration $t$ is some $Y_t$ modified by the scaling matrix. Vector $Y_t$ is drawn from a family of random variables $Y$ defined for all $\eta$, and $Y_t=Y(\eta_t)$. The result is stated here.

\begin{theorem}[Theorem 3.2 in \cite{sunehag_variable_2009}]
\label{thm:sunehag}
Let $\mathbb{L}:\R^k\rightarrow \R$ be a twice differentiable function with a unique minimum  $\overline\eta$ and
$\eta_{t+1} = \eta_t - \gamma_t B_t Y_t$ where $B_t$ is symmetric and only depends on information available at time $t$. Then $\eta_t$ converges to $\overline\eta$ almost surely if the following conditions hold
\begin{equation}\nonumber
\begin{split}
\textbf{C.1}\hspace{0.1cm} &(\forall t)\hspace{0.1cm} \E_t Y_t = \nabla \mathbb{L}(\eta_t)\\
\textbf{C.2}\hspace{0.1cm} &(\exists K) (\forall \eta)\hspace{0.1cm} \norm{\nabla^2_\eta \mathbb{L}(\eta)}\leq 2K \\
\textbf{C.3}\hspace{0.1cm} &(\forall\delta>0)\hspace{0.1cm} \inf_{ \mathbb{L}(\eta)-\mathbb{L}(\overline \eta)>\delta}\norm{\nabla \mathbb{L}(\eta)}>0\\
\textbf{C.4}\hspace{0.1cm} &(\exists A,B)(\forall t)\hspace{0.1cm}\E_t\norm{ Y_t}^2
\leq A+B\cdot \mathbb{L}(\eta_t)\\
\textbf{C.5}\hspace{0.1cm} &(\exists m,M:0 < m < M < \infty )\hspace{0.1cm} (\forall t) mI \prec B_t \prec MI, \text{where } I \text{ is the identity matrix;} \\ 
\textbf{C.6} \hspace{0.1cm}&  \sum_t \gamma_t^2 <\infty, \hspace{0.1cm}\sum_t \gamma_t =\infty
\end{split}
\end{equation}
\end{theorem}
$\E_t$ in conditions \textbf{C.1} and \textbf{C.4} notes the conditional expectation given observations until time $t$. That is, $\E_t X=\E\left[ X\mid \F_t\right]$, where $\F_t =\{\eta_1,...,\eta_t\}$ in this case. We recall Robbins-Siegmund Theorem \cite{robbins_siegmund_1971} below, which is the key tool for proving both Theorem~\ref{thm:sunehag} and our generalization.

\begin{theorem}[Robbins-Siegmund]\label{theorem:Robbins-Siegmund}
Let $(\Omega, \F, P)$ be a probability space and $\F_1\subseteq \F_2\subseteq \cdots$ a sequence of sub-$\sigma$-fields of $\F$. Let $U_t, \beta_t, \epsilon_t$ and $\zeta_t$, $t=1,2,...$ be non-negative $\F_t$-measurable random variables such that 
\begin{align}\label{eq:Robbins-Siegmund}
    \E(U_{t+1}\mid \F_t)\leq (1+\beta_t)U_t + \epsilon_t-\zeta_t, \ \ t=1,2,...
\end{align}
Then on the set $\left\{\sum_t \beta_t<\infty, \sum_t \epsilon_t<\infty\right\}$, $U_t$ converges almost surely to a random variable, and $\sum_t\zeta_t<\infty$ almost surely.
\end{theorem}

\section{Dual Stochastic Natural Gradient Descent}\label{sec:dsngd}


Recall from the introduction that the main idea to reduce the computational complexity of DSNGD is to define our optimization problem over a DFM. We start by establishing in section~\ref{sec:joint} that the family of joint distributions $P(\Y,\X)$ satisfying the core MLR assumption is a LEF and hence a DFM. Then, we rely on duality to provide an efficient computation of the natural gradient of the log-loss function in section~\ref{sec:fast-natural-gradient}. Finally, we provide the DSNGD algorithm in section~\ref{sec:dsngd_definition}. 


\subsection{MLR generative model. The joint distribution\label{sec:joint}}

The next result proves that the the family of joint distributions $P(\Y,\X)$ satisfying the core MLR assumption is a LEF and hence a DFM.
\begin{proposition}\label{propo:model}
The log-odds ratio of the class posteriors $P(\Y\mid\X)$ is an affine function of the features $\X$ if and only if the joint distribution $P(\Y,\X)$ belongs to LEF.

Furthermore, there exist the LEF natural parametrization of the joint distribution
 \begin{equation}\label{eq:natural_parametriztion}
\begin{split}
    P_{\eta}(x,y)&=\frac{\exp{S(y)^\intercal\alpha+T(x)^\intercal\beta_y }}{\lambda(\eta)}\\
    \lambda(\eta)& = \int_x \sum_y \exp{S(y)^\intercal\alpha+T(x)^\intercal\beta_y }
\end{split}
\end{equation}
where $\eta=(\alpha,\beta)$,  $\alpha\in\R^{s-1}$, $\beta\in\R^{s \times t}$, $\beta_y$ is the $y$-th row of $\beta$ and
\begin{equation}
    \begin{split}
    T:\Omega&\rightarrow \R^{t}\\
    S:[1,...s]&\rightarrow \R^{s-1}
    \end{split}
\end{equation} are sufficient and minimal statistics of $\X$ and $\Y$ respectively.
\end{proposition}

The proof relies strongly on theorem 2 in \cite{banerjee_analysis_2007} and can be found in appendix~$\ref{proof:lexyf_parametrization}$.
This is convenient for our purpose, because, if we recall section \ref{dfm}, in a DFM the costs of natural gradient computations can be highly reduced, based on the property shown by equation \ref{eq:natural_gradient_is_gradient}. Next, we provide the dually flat parametrization of $P(\Y,\X)$.

\subsubsection{Dually flat parametrization of the joint distribution}

We have seen that $P(\Y,\X)$ is a LEF and that we can choose the natural parametrization of equation~\ref{eq:natural_parametriztion}. With a linear transformation, $S$ can become a canonical statistic,  that means, $S(i)_j=\delta_{i=j}$ for $1\leq i<s$ and $S(s)=0$. For simplicity, we fix statistic $S$ to be canonical from now on. The conditional probability distributions with $\eta$ parametrization are
\begin{equation}\label{eq:conditionals}
    P_{\eta}(y\mid x)=\frac{\exp{S(y)^\intercal\alpha+T(x)^\intercal\beta_y }}{\sum_y \exp{S(y)^\intercal\alpha+T(x)^\intercal\beta_y }}
\end{equation}

As \cite{amari_information_2016} proves, the exponential family manifold is built after the convex function $F(\eta)=\log \lambda(\eta)$. The reference proves that this Riemannian manifold derived from  $F(\eta)$, according to equation \ref{eq:F_manifold}, has the Fisher information metric, as is usually considered for statistical manifolds, defined as
\begin{equation}
    G_\eta= -\E_{x,y\sim P_\eta} \left[ \nabla^2 \log P_\eta (y,x)  \right]
\end{equation}


Equation~\ref{eq:F_manifold} also reveals the dual parametrization $\eta^*$. For LEF, it is called the expectation parametrization and it is shown below. For more properties of the dual parametrization see \cite{amari_information_2016}.
To simplify the notation, if $x=\begin{pmatrix}
x_1 &\cdots &x_n 
\end{pmatrix} $, we note $\nabla_x=
\begin{pmatrix}
\frac{\partial}{\partial x_1} &\cdots &\frac{\partial}{\partial x_n} 
\end{pmatrix}^\intercal
$. So for every $i\in\{1,...,s\}$ write

\begin{equation}\label{eq:expect_parametriztion}
\begin{split}
        \alpha^* =\nabla_\alpha F(\eta)& =\sum_y S(y)P_\eta(y)=  \E_{\Y}[S(y)]=(P_\eta(\Y=1),...,(P_\eta(\Y=s-1))^\intercal\\
    \beta_{i}^* =\nabla_{\beta_{i}} F(\eta)& =P_\eta(\Y=i)\int_\X T(x)P_\eta(x\mid \Y=i)=  P_\eta(\Y=i)\E_{\X\mid \Y=i}[T(x)]
\end{split}
\end{equation}
Define $\eta^*=(\alpha^*,\beta^*)$ with $\beta = (\beta_1^*,...,\beta_s^*)$ the dual parameterization, or equivalently, the expectation parameters.

Observe that $P(\Y)$ is the categorical distribution (since $\Y$ is  discrete and finite) and therefore it is a LEF, where $\alpha^*$ are actually the expectation parameters. Moreover
\begin{align}\label{eq:expect_param_x|y}
    \theta_i\coloneqq\theta_i(\alpha^*,\beta^*_i)= \frac{\beta^*_{i}}{P_{\eta^*}(\Y=i)}= \E_{\X\mid \Y=i}[T(x)]
\end{align} are the expectation parameters of the conditional distribution $P(\X\mid\Y = i)$.

\subsection{Fast natural gradient of the log-loss\label{sec:fast-natural-gradient}}
  This section allows to compute the natural gradient of the log-loss function without having to use the metric matrix directly, but using both dual parametrizations instead. 
  
Given $(y,x)\in\Y\times\X$ and $\eta\in\R^k$, the log-loss function is defined as
 \begin{equation}
l(\eta,x,y) = -\log P_\eta(y\mid x)
\end{equation}

Below result reveals $\widetilde\nabla l(\eta,x,y)$ using both dual parametrizations $\eta$ and $\eta^*$.

\begin{proposition}\label{propo:nat_grad_dsngd}
Let $l$ be the log-loss function. Then, if $P(\Y,\X)$ is a
DFM 
\begin{equation}
    \label{eq:nat_grad_dsngd}
    \widetilde{\nabla}l(\eta,x,y) = \nabla h(x,\eta^*)\cdot (q_\Y(x,P_{\eta})-e_s(y))
\end{equation}
where
\begin{equation}
q_\Y(x,P)=\begin{pmatrix}P(\Y=1|x)\\
\vdots\\
P(\Y=s|x)
\end{pmatrix},
\end{equation}
$h(x,\eta^*)=(\log P_{\eta^*}(\Y=1,x),...\log P_{\eta^*}(\Y=s,x))$ and $e_s(k)$ is the $k$-th canonical $s$-dimensional vector.
\end{proposition}

The proof of proposition~\ref{propo:nat_grad_dsngd} is presented in Appendix \ref{proof:nat_grad_dsngd}. To evaluate the computational complexity of using equation \ref{eq:nat_grad_dsngd} we determine an expression of $\nabla h(x,\zeta^*)$ with respect to the expectation parameters $\theta_y$ of $\X$ given $\Y$ already mentioned in equation~\ref{eq:expect_param_x|y}. Below notation is used
\begin{equation}
K_i=\left(
    \begin{array}{ccc|c}
    & &  &-1 \\
    & Id^{i-1} & & \vdots\\
    & &  & -1\\
    \end{array}\right)\hspace{1cm}
    d(x,y,\zeta^*)=\frac{1-\theta_y^\intercal \nabla_{\theta_y}\log P_{\theta_y}(x\mid y)}{P_{\zeta^*}(y)}
\end{equation} and the proof is shown in Appendix \ref{proof:h}.

\begin{proposition}\label{theorem:h}
\begin{equation}\nabla h(x,\zeta^*)=
    \begin{pmatrix}
     \nabla_{\alpha^*} h(x,\zeta^*)\\
     \nabla_{\beta_1^*} h(x,\zeta^*)  \\
     \vdots\\
     \nabla_{\beta_s^*} h(x,\zeta^*)  
    \end{pmatrix}
\end{equation} where
\begin{equation}
\label{eq:generic_nabla_h}
\begin{split}
    \nabla_{\alpha^*}h(x,\zeta^*)=& K_s\cdot diag(d(x,1,\zeta^*),...,d(x,s,\zeta^*))
    \\
    \nabla_{\beta_k^*}h(x,\zeta^*)=&\frac{\nabla_{\theta_k}\log P_{\theta_k}(x\mid \Y=k)\cdot e_s(k)^\intercal}{P_{\zeta^*}(\Y=k)}
\end{split}
\end{equation}
\end{proposition}

The complexity analysis of natural gradient is presented now, and the reader can find the proof in appendix~\ref{proof:natural_gradient_complexity}.

\begin{proposition}\label{propo:nat_grad_complexity}
The computational complexity of the natural gradient $\widetilde{\nabla}l(\eta,x,y)$ using proposition~\ref{propo:nat_grad_dsngd} is $O(s\cdot( A+t))$ where $A$ is the cost of computing $\nabla_{\theta_y}\log P_{\theta_y}(x\mid y)$, $s$ is the number of classes and $t$ is the dimension of statistic $T$. 
\end{proposition}

Observe that the manifold dimension is $k=s-1+s\cdot t$ and therefore, a computation is linear on the number of the variables of the model if its complexity order is $O(k) =O(s(1+t))=O(st)$. Therefore, the costs of computing the natural gradient can be reduced to linear if the cost $A$ is low enough, precisely, if $A$ is at most linear ($O(A)\leq O(k)$). This is the case when $\X$ is discrete and finite (section \ref{sec:dsngd_complexity}). 


\subsection{DSNGD definition}\label{sec:dsngd_definition}
DSNGD aims to solve the MLR optimization problem using the natural parametrization $\eta$ of the LEF on $\Y\times\X$: If $\overline P$ is an unknown probability distribution over $\Y\times\X$, optimize $\mathbb{L}(\eta)=\E_{x,y\sim \overline P}\left[l(\eta,x,y)\right]$ for $\eta\in\R^k $ where $l(\eta,x,y)$ is the log-loss function. 
The solution $\overline\eta\in\R^k$ to this problem refers to the conditional distributions $P_{\overline{\eta}}(\Y\mid \X)$ that better fits the hidden conditional distributions $\overline P(\Y|\X) $. To that end, we define a stochastic natural gradient based algorithm.

Using proposition $\ref{propo:nat_grad_dsngd}$, DSNGD moves by following the update equation
\begin{definition}[DSNGD update]
\begin{equation}\label{eq:dsngd_update}
\eta_{t+1}=\eta_{t}-\gamma_{t}\nabla h(x_t,\zeta_t^*)\cdot (q_\Y(x_t,P_{\eta_t})-e_s(y_t))
\end{equation}
where $\left\{ \zeta_{t}^{*}\right\} _{t\in\mathbb{N}}$ is a sequence in the expectation parametrization such that $\left\{ \zeta_{t}\right\} _{t\in\mathbb{N}}$ converges.
\end{definition}
Note that $q_\Y(x_t,P_{\eta_t})$ is a stable term (it only takes values between $0$ and $1$). Moreover, DSNGD forces the stability of the $\nabla h(x_t,\zeta_t^*)$ term, since $\zeta_t$ is a convergent sequence. This is the same strategy of CSNGD, and similarly, it is going to ensure the convergence of the algorithm in theorem~\ref{theorem:dsngd_convergence}. Observe that equation~\ref{eq:dsngd_update} is also well defined when the parameterization is not minimal (when $T$ is not a minimal statistic), therefore DSNGD can be run in such general case, where $S$ and $T$ are not minimal. Steps taken by DSNGD are specified in Algorithm \ref{alg:DSNGD} below.

\begin{algorithm}[H]\label{alg:DSNGD}
\SetAlgoLined
\KwResult{$\eta$ }
 $\eta\leftarrow\eta_0,\zeta^*\leftarrow\zeta_0^*,\gamma\leftarrow\gamma_0$\;
 \While{observations $x,y$ and stopping condition is false}{
  $q\leftarrow q_\Y(x,P_{\eta})$\;
  $grad\_h\leftarrow \nabla h(x,\zeta^*)$\;
  $d\leftarrow grad\_h\cdot (q-e_s(y))$\;
  $\eta \leftarrow \eta -\gamma \cdot d$\;
  update $\zeta^*$\;
  update $\gamma$
 }
 \caption{DSNGD}
\end{algorithm}

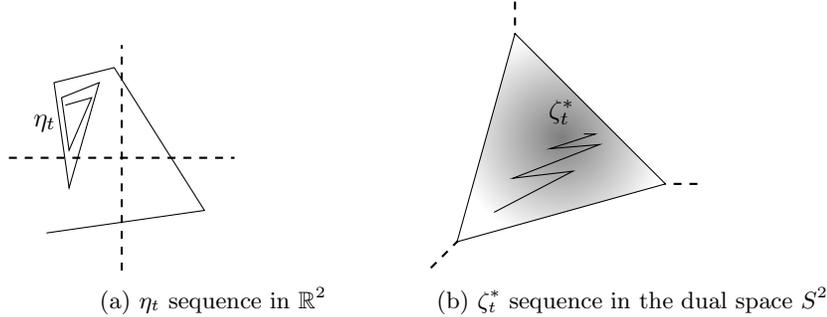
\begin{figure}[h] 
    \begin{subfigure}[b]{0.45\textwidth}
        \begin{tikzpicture}
            \draw[thick, dashed] (-1.5,0) -- (1.5,0);
            \draw[thick, dashed] (0,-1.5) -- (0,1.5);
            \draw (-1,-1) -- (1.1, -0.7) --(-0.1,1.2)--(-0.9,1)--(-0.7,-0.4)--(-0.3,1)--(-0.8,0.8)--(-0.7,0.1)--(-0.4,0.8)--(-0.75,0.7) node[anchor=north east] {$\eta_t$};
        \end{tikzpicture}
         \caption{$\eta_t$ sequence in $\R^2$}
         \label{fig:R2}
     \end{subfigure}
     \begin{subfigure}[b]{0.45\textwidth}
        \begin{tikzpicture}
            \draw[thick, dashed] (0,0,0) -- (2.5,0,0);
            \draw[thick, dashed] (0,0,0) -- (0,2.5,0);
            \draw[thick, dashed] (0,0,0) -- (0,0,3);
            \shadedraw[shading=radial, color=lightgray, draw=black] (2,0,0) -- (0,2,0)--(0,0,2)--cycle;
            \draw (0.3,0.2,1.5)-- (1.,0.4,0.6)--(0.4,0.5,1.1)--(1.2,0.6,0.2)--(0.7,0.7,0.6) --(1.1,0.7,0.1)  --(1.0,0.75,0.22)node[anchor=south east] {$\zeta^*_t$};
        \end{tikzpicture}
         \caption{$\zeta^*_t$ sequence in the dual space $S^2$}
         \label{fig:S2}
    \end{subfigure}
    \caption{$\eta_t$ and $\zeta^*_t$ sequences obtained in DSNGD where $\X=\{0\}$ and $\Y=\{0,1,2\}$
    }
            \label{fig:1}
\end{figure}

The sequence $\left\{ \zeta_{t}^{*}\right\} _{t\in\mathbb{N}}$, or simply $\zeta^*_t$ as an abuse of notation, can be any sequence in the dual space whose dualized sequence $\left\{ \zeta_{t}\right\} _{t\in\mathbb{N}}$ is convergent. For example, it can be constant. The resulting algorithm keeps track of two independent sequences; the main sequence $\eta_{t}$ which estimates the solution $\overline{\eta}$ to the problem, and the sequence $\zeta_{t}^{*}$ selected with the convergence constraint and whose space is the dual. For example, assume the trivial case where $\X=\{0\}$ and $\Y=\{0,1,2\}$. The only conditional probability distribution of the problem is the Categorical distribution $P(\Y\mid \X=0)$. This space is represented by $\R^2$ and its dual space is represented by the simplex $S^2$. Then, the main sequence $\eta_t$ moves in $\R^2$ while the independent sequence $\zeta^*_t$ traces its path in $S^2$. Figure $\ref{fig:1}$ illustrates iterations followed by $\eta_t$ (instruction line $6$ of the algorithm) and $\zeta_{t}^{*}$ (instruction line $7$ of the algorithm) when running DSNGD for this simple example.

Recall that the sequence $\zeta_t^*$ can be chosen freely as long as its dual is convergent. However, recall that DSNGD is a natural gradient based algorithm. The algorithm effectively takes a natural gradient step only when $\eta_t$ and $\zeta^*_t$ refer to the same probability distribution point, according to equation \ref{eq:dsngd_update} and proposition $\ref{propo:nat_grad_dsngd}$. In section \ref{sec:dsngd_convergence} there is our proof of DSNGD convergence to the solution $\overline{\eta}$, and if $\zeta_t^*$ is selected such that it also converges to the solution, then both sequences get closer along the optimization process, turning DSNGD steps into more accurate approximations of natural gradient steps. Therefore, in order to benefit from natural gradient speed up properties, it is recommended that sequence $\zeta_t^*$ converges to the solution $\overline{\eta}^*=\nabla F(\overline{\eta})$. For example, this can be accomplished by determining $\zeta^*_t$ using a maximum a posteriori estimator of the parameters of $P(\Y,\X)$ obtained from data up to $t$.  


\section{Discrete DSNGD and computational complexity}\label{sec:dsngd_complexity}

This section assumes that space $\X$ is discrete, that is $\X=\{1,...,m\}$ for some $m\in\N$. For simplicity, we assume  $T$ to be also a canonical statistic, that is, $T(i)_j=\delta_{i=j}\in\R^{m-1}$ for $1\leq i<m$ and $T(m)=0$. A theorem deduces and proves that the complexity order for discrete DSNGD of one iteration is linear on the dimension of the parameter $\eta$. Let us show a simple example of discrete DSNGD to begin with.

\subsection{Example}
Let $\Y=\{1,2\}$ and $\X=\{1,2\}$ and minimal and canonical statistics $S$ and $T$. Let  $\eta=(\alpha,\beta)$ be the natural parameter and $\zeta^*=(\alpha^*, \beta^*)$ be the independent dual parameter. Observe that in this case, $\alpha$ and $\alpha^*$ are $1$-element vectors and $\beta$ and $\beta^*$ are $2$-element not squared matrices. In this example we complete an iteration of discrete DSNGD algorithm, following the instructions listed in algorithm $\ref{alg:DSNGD}$.

Let $(y,x)=(2,1)$ be an observation. Statistics $T$ and $S$ are assumed to be canonical. Instruction line $3$ consist on using equation $\ref{eq:conditionals}$ to compute
\begin{equation}
    q_\Y(x,P_\eta) =
    \begin{pmatrix}
    P(\Y=1\mid x)\\
    P(\Y=2\mid x)
    \end{pmatrix} =R\cdot
    \begin{pmatrix}
    \exp{(\alpha_1+\beta_1)}\\
    \exp{\beta_2}
    \end{pmatrix}
\end{equation}
where $R=\frac{1}{\exp{(\alpha_1+\beta_1)}+\exp{\beta_2}}$.
For instruction line $4$, express function $h(x,\zeta^*)$ (use equation~\ref{eq:expect_parametriztion}), then apply the gradient.
\begin{equation}
    h(x=1,\zeta^*) =(\log \beta_1^*, \log \beta_2^*)\rightarrow \nabla h(x=1,\zeta^*)=
    \begin{pmatrix}
    0 & 0\\
    \frac{1}{\beta_1^*} & 0\\
    0 & \frac{1}{\beta_2^*}
    \end{pmatrix}
\end{equation}
Proceed now with instruction line $5$. It computes the approximation of the natural gradient and the direction that DSNGD uses for the $\eta$ update.
\begin{equation}
\begin{split}
    d=&  \nabla h(x=1,\zeta^*) (q_\Y(x,P_\eta) - e_2(y=2))\\
    =&
    \begin{pmatrix}
    0 & 0\\
    \frac{1}{\beta_1^*} & 0\\
    0 & \frac{1}{\beta_2^*}
    \end{pmatrix}\cdot
        \begin{pmatrix}
    R\cdot\exp{(\alpha_1+\beta_1)}\\
    (R\cdot\exp{\beta_2}) -1
    \end{pmatrix}\\
    =& 
    \begin{pmatrix}
    0 \\
    \frac{R\cdot\exp{(\alpha_1+\beta_1)}}{\beta_1^*}\\
    \frac{(R\cdot\exp{\beta_2}) -1}{\beta_2^*}
    \end{pmatrix}
\end{split}
\end{equation}
Next instruction lines of the algorithm are standard to update the parameter vector $(\alpha,\beta_1,\beta_2)$ using direction $d$, so there is no need to go further. If for instance the observation is $(y,x)=(2,2)$, then the approximation of the natural gradient is
\begin{equation}
    d=  
    \begin{pmatrix}
    \frac{R\exp\alpha_1}{\alpha_1^*-\beta_1^*} -\frac{R-1}{1-\alpha_1^*-\beta_2^*} \\
    \frac{-R\cdot\exp{\alpha_1}}{\alpha_1^*-\beta_1^*}\\
    \frac{R -1}{1-\alpha^*_1-\beta_2^*}
    \end{pmatrix}
\end{equation}where $R=\frac{1}{1+\exp{\alpha_1}}$






Before analyzing the computational complexity of DSNGD, it is necessary to determine the generator of $\zeta^*_t$ sequence. Sequence $\zeta^*_t$ belongs to the dual space of the LEF distributions on $\Y\times\X$, and if $S$ and $T$ are canonical statistics then it implies that $\zeta^*_t$ are directly the probabilities $P(y,x)$ after equation $\ref{eq:expect_parametriztion}$. It is possible to select the well known maximum a posteriori (MAP) estimator with parameter $a\in\R$. This estimator is a simple counting of observations over the discrete space $\Y\times\X$ with an starting assumption of incidence of $a$ for every event $y,x$. This estimator is linear and it clearly converges (to the solution). 

 First a similar result as proposition~\ref{theorem:h} is stated, taking into account the new assumption on $\X$. The proof of proposition below is found in appendix \ref{proof:h_discrete}.

\begin{proposition}\label{theorem:h_discrete}Let $\X=\{1,...,m\}$ and let $T$  be a minimal and canonical statistic. Then
\begin{equation}
    \label{eq:discrete_nabla_h}
    \begin{split}
        \nabla_{\alpha^*} h(x,\zeta^*)&=
        \begin{cases}0 & x\neq m \\
    K_s\cdot diag(\frac{1}{
        P_{\zeta^*}(x,\Y= 1)},...,\frac{1}{
        P_{\zeta^*}(x,\Y= s)}) & x=m
        \end{cases}\\
        \nabla_{\beta_y^*}h(x,\zeta^*)&=\frac{1}{
        P_{\zeta^*}(x, y)}\cdot
        \begin{cases}
        e_{m-1}(x)\cdot e_s(y)^\intercal & x\neq m\\
        -\textbf{1}_{m-1}
        \cdot e_s(y)^\intercal & x= m
        \end{cases}
    \end{split}
\end{equation}where $\textbf{1}_n\in\R^n$ is a vector filled with ones at every coordinate.  
\end{proposition}

Now it is possible to analyze the computational complexity of discrete DSNGD. Below theorem proves that DSNGD, just as SGD, is a linear algorithm.
\begin{theorem}
Let $\X=\{1,...,m\}$ and let $T$  be a minimal and canonical statistic. Assume estimator $\zeta^*$ of DSNGD is linear. Then discrete DSNGD iterations have linear complexity order on the manifold dimension.
\end{theorem}

\begin{proof}
Let $k=(s-1) + s\cdot t$ be the dimension of $\eta$. Then $O(k) = O(st)$. Analyze the computational complexity of discrete DSNGD. That is, analyze the computational cost of instruction lines $3,4,5,6$ and $7$ shown in Algorithm $\ref{alg:DSNGD}$. 

Complexity of instruction lines $3,4$ and $5$ is given by proposition~\ref{propo:nat_grad_complexity}, which is $O(sA+st))$ where $s A$ is the cost of computing $\nabla_{\theta_y}\log P_{\theta_y}(x\mid y)$ for all $y\in\Y$. Observe equations~\ref{eq:expect_parametriztion} and \ref{eq:expect_param_x|y} assuming $T$ canonical and write
\begin{equation}
\begin{split}
        \alpha^* &=(P_{\zeta^*}(\Y=1),...,P_{\zeta^*}(\Y=s-1))^\intercal\\
    \theta_y&=(P_{\zeta^*}(\X=1\mid y),...,P_{\zeta^*}(\X=m-1\mid y))^\intercal
\end{split}
\end{equation}Deduce then that $O(sA)= O(k)$.


Instruction line $6$ adds $k$ operations.

Finally, recall that a linear complexity order estimator is chosen for $\zeta^*_t$ sequence, implying that instruction line $7$ is of linear order $O(k)$.

In conclusion, the computational complexity order of DSNGD is \begin{equation}
     O(sA +st) +O(k)+ O(k)= O(k)
\end{equation} and therefore linear. 
\end{proof}
\section{Discrete DSNGD convergence}\label{sec:dsngd_convergence}
In this section we prove the convergence of the discrete DSNGD. Discrete DSNGD refers to the case where $\X=\{1,...,m\}$ for some $m\in\R$. We start by generalizing Theorem~3.2 in Sunehag's et al. \cite{sunehag_variable_2009} (introduced above and referred to from now on as Theorem~\ref{thm:sunehag}) in Section~\ref{sec:convergence_generalization}. This generalization provides enough flexibility so as to be used later to prove the convergence of DSNGD in Section~\ref{sec:dsngd_convergence_2}. 
\subsection{Generalizing Sunehag et. al. variable metric stochastic approximation theory.}\label{sec:convergence_generalization}

Theorem \ref{thm:sunehag} is used to prove CSNGD convergence, however it can not be used to prove DSNGD convergence. First, because it demands the vector it follows to be factored as the product of a matrix $B_t$ and a vector $Y_t$ that approximates the gradient (condition \textbf{C.1}). But DSNGD is defined to directly approximate the natural gradient, without the gradient as reference.  And second, even if DSNGD is written as the product of a matrix and a vector, matrix $\nabla h(x_t,\zeta_t^*)$ is not squared. So we need a more general convergence theorem. 

Our result proves almost sure convergence of the sequence
\begin{equation}
    \eta_{t+1} = \eta_t - \gamma_t  Y(\eta_t,\F_t)
\end{equation} where $Y(\eta,\F_t)$ is a family of random vectors defined for every $\eta$ and for every set \begin{equation}\label{eq:F}
    \F_t=\{(y_i,x_i)\mid i<t\}
\end{equation}As an abuse of notation write $Y_t$ meaning the random variable $Y(\eta_t,\F_t)\in\R^n$. The main modification with respect to Theorem~\ref{thm:sunehag} is that we unify conditions C.1 and C.3
\begin{equation}
\begin{split}
\textbf{C.1}\hspace{0.5cm} & (\forall t)\hspace{0.5cm} \E_t Y_t = \nabla l(\eta_t)\\
\textbf{C.3}\hspace{0.5cm} & (\forall\delta>0)\hspace{0.5cm} \inf_{ l(\eta)-l(\overline \eta)>\delta}\norm{\nabla l(\eta)}>0
\end{split}
\end{equation}
to instead require 
\begin{equation}\textbf{C.3}\hspace{0.5cm} (\forall\delta>0)\hspace{0.5cm} \inf_{ l(\eta_t)-l(\overline \eta)>\delta}\nabla l(\eta_t)^T\E_t \left[ Y_t\right]>0
\end{equation}
New condition \textbf{C.3} uses $\E_t$, referring to the conditional expectation given $\F_t$ of equation~\ref{eq:F}, which is a generalization of the definition of $\E_t$ in Sunehag \cite{sunehag_variable_2009}.

Theorem \ref{thm:sunehag} imposes that the expectation of the step taken must be the gradient and that the norm of the gradient must not 
approach to zero outside any
environment of the minimum. Instead, we impose that the expectation of the step taken must not
approach to the border
of the half-space which has the gradient as its normal vector, unless we are approaching the minimum simultaneously. This is a more general condition.
Furthermore, condition \textbf{C.5} on the maximum and minimum eigenvalues of the matrix $B_t$ can also be removed. In fact, our result proves the convergence of algorithms with scaling matrices $B_t$ whose spectrum is not bounded from below by a strictly positive number, as long as new version of condition \textbf{C.3} holds.

It is formally stated below. Proof is found in appendix \ref{proof:convergence_generalization}.  
\begin{theorem}\label{theorem:convergence}
Let $l:\R^k\rightarrow \R$ be a twice differentiable function with a unique minimum  $\overline\eta$ and
$\eta_{t+1} = \eta_t - \gamma_t  Y_t$. Then $\eta_t$ converges to $\overline\eta$ almost surely if the following conditions hold
\begin{equation}\nonumber
\begin{split}
\textbf{C.2}\hspace{0.1cm} &(\exists K) (\forall \eta)\hspace{0.1cm} \|\nabla^2_\eta l(\eta)\|\leq 2K \\
\textbf{C.3}\hspace{0.1cm} &(\forall\delta>0)\hspace{0.1cm} \inf_{ l(\eta_t)-l(\overline\eta)>\delta}\nabla l(\eta_t)^T\E \left[ Y_t\right]>0\\
\textbf{C.4}\hspace{0.1cm} &(\exists A,B)(\forall t)\hspace{0.1cm}\E\| Y_t\|^2
\leq A+Bl(\eta_t)\\
\textbf{C.6} \hspace{0.1cm}&  \sum_t (\gamma_t)^2 <\infty, \hspace{0.1cm}\sum_t \gamma_t =\infty
\end{split}
\end{equation}
\end{theorem}

\subsection{Proving convergence}\label{sec:dsngd_convergence_2}

Next, we show how Theorem~\ref{theorem:convergence} can be used to prove DSNGD convergence in the discrete case. That is, we use it to prove the next result:

\begin{theorem}\label{theorem:dsngd_convergence}
DSNGD in the canonical parametrization such that $\Y$ and $\X$ are discrete, converges almost surely to the optimum.
\end{theorem}

The proof consists on showing that conditions \textbf{C.2, C.3, C.4} and \textbf{C.6} of theorem~\ref{theorem:convergence} hold. Condition \textbf{C.6} is assumed to hold, by just selecting an appropriate sequence of learning rates $\gamma_t$. Conditions \textbf{C.2} and \textbf{C.4} are proved in appendices \ref{proof:C2} and \ref{proof:C4} respectively. Proof of condition \textbf{C.3} is shown below.
\begin{proof}

Compute the gradient of $l(\eta)$ (see equation~\ref{eq:log-loss_grad}) and use proposition \ref{propo:nat_grad_dsngd} to obtain  $\E_t \left[Y_t\right]$ involved in condition \textbf{C.3}.

\begin{equation}\label{eq:1}
    \begin{split}
        \nabla l(\eta) = & \E_t \left[\nabla l(\eta,x,y)\right] \\
        = & \sum_x \nabla h(x,\eta)\sum_y (q_\Y(x)-e_s(y))\overline P(x,y)\\
        = & \sum_x \nabla h(x,\eta) diff_\Y(x,\eta)\\
         \E_t \left[Y_t\right] =& \sum_x \nabla h(x,\zeta^*) diff_\Y(x,\eta)
    \end{split}
\end{equation}
where 
\begin{equation}
    diff_\Y(x,\eta)= (q_\Y(x,P_\eta)-q_\Y(x,\overline P))\overline P(x)
\end{equation}

Further evolve equation \ref{eq:1} to finally multiply $\nabla l(\eta)^\intercal\E_t \left[Y_t\right]$ and check condition~\textbf{C.3}. Continue by developing $\nabla l(\eta)$ first, precisely compute $\nabla h(x,\eta)$. To simplify the notation, decompose $\nabla = (\nabla_\alpha,\nabla_{\beta_1},...,\nabla_{\beta_s})$

\begin{equation}
    \begin{split}
        \nabla_{\alpha} h(x,\eta)&=S+u(P_\eta)\cdot (1,...,1)\hspace{2cm} u(P) = -\sum_{y} S(y)P(y)\\
        \nabla_{\beta_y} h(x,\eta)&= T(x) e_s(y)^\intercal+v(y,P_\eta)\cdot (1,...,1)\hspace{0,5cm} v(y,P) = -\sum_{x}T(x)P(y,x)
    \end{split}
\end{equation}

Since $(1,...,1)\cdot diff_\Y(x,\eta)=0$ then

\begin{equation}
    \begin{split}
         \nabla_{\alpha} l(\eta)= & \sum_x \nabla_{\alpha} h(x,\eta) diff_\Y(x,\eta)\\
         = & \sum_x  S\cdot diff_\Y(x,\eta)\\
         = & S \cdot diff_\Y(\eta)\\
         \nabla_{\beta_y} l(\eta)= & \sum_x \nabla_{\beta_y} h(x,\eta) diff_\Y(x,\eta)\\
         = & \sum_x  T(x) e_s(y)^\intercal diff_\Y(x,\eta)\\
         =& T\cdot diff_\X(y,\eta)
    \end{split}
\end{equation}

where
\begin{equation}
    \begin{split}
        diff_\Y(\eta) = &  
        \begin{pmatrix}
        P_\eta(\Y=1)-\overline P(\Y=1)\\
        \vdots\\
        P_\eta(\Y=s)-\overline P(\Y=s)
        \end{pmatrix}
        \\
        diff_\X(y,\eta) = &
       \begin{pmatrix}
        (P_\eta(\Y=y\mid\X=1)-\overline P(\Y=y\mid\X=1))\overline P(\X=1)\\
        \vdots\\
        (P_\eta(\Y=y,\X=m)-\overline P(\Y=y\mid \X=m))\overline P(\X=m)
        \end{pmatrix}
    \end{split}
\end{equation}
 
Now develop $\E_t Y_t $ further. Recall that the canonical parametrization is selected so plug in  proposition \ref{theorem:h_discrete} into equation $\ref{eq:1}$. Decompose $ \E_t = ( \E_{t,\alpha^*}, \E_{t,\beta_1^*},..., \E_{t,\beta_s^*})$
\begin{equation}
    \begin{split}
         \E_{t,\alpha^*} \left[Y_t\right]= & \sum_x \nabla_{\alpha^*} h(x,\zeta^*) diff_\Y(x,\eta)\\
         = & K_{s}\cdot diag(d(m,1,\zeta^*),...,d(m,s,\zeta^*))
        \cdot diff_\Y(m,\eta)\\
         \E_{t,\beta^*_y} \left[Y_t\right]= & \sum_x \nabla_{\beta_y^*} h(x,\zeta^*) diff_\Y(x,\eta)\\
         =&K_{m}\cdot diag(
        d(1,y,\zeta^*),...,d(m,y,\zeta^*))\cdot diff_\X(y,\eta)
    \end{split}
\end{equation}

Proceed now to check the condition. Develop the products until obtain

\begin{equation}
\begin{split}
        \nabla_\alpha l(\eta)^\intercal\E_{t,\alpha^*} \left[Y_t\right]=&\sum_y c(y)\\
        \nabla_{\beta_y} l(\eta)^\intercal\E_{t,\beta_y^*} \left[Y_t\right]
        =&-c(y)+\sum_x d(x,y,\zeta^*)(P_\eta(y|x)-\overline P(y|x))^2 \overline{P}(x)^2
    \end{split}
\end{equation}
where $c(y) = d(m,y,\zeta^*)(P_\eta(y)-\overline P(y))(P_\eta(y|x=m)-\overline P(y|x=m))\overline P(x=m)$

Finally,
\begin{equation}\label{eq:c4}
    \begin{split}
                \nabla l(\eta)^\intercal\E_t \left[Y_t\right]=&\nabla_\alpha l(\eta)^\intercal\E_{t,\alpha^*} \left[Y_t\right]+\sum_y \nabla_{\beta_y} l(\eta)^\intercal\E_{t,\beta_y^*} \left[Y_t\right]\\
                =&\sum_y c(y)+\sum_y-c(y)+\sum_x d(x,y,\zeta^*)(P_\eta(y|x)-\overline P(y|x))^2\overline P(x)^2\\
                =& \sum_{y,x} d(x,y,\zeta^*)(P_\eta(y|x)-\overline P(y|x))^2\overline P(x)^2
    \end{split}
\end{equation}
Notice in equation \ref{eq:c4} that $\nabla l(\eta)^\intercal\E_t \left[Y_t\right]$ is a sum of positive numbers, and it vanishes only if $\eta=\overline\eta$. Also, since $d(x,y,\zeta^*)> 1$, observe that 

\begin{equation}
\begin{split}\label{eq:c3_proof_key}
     \nabla l(\eta)^\intercal\E_t \left[Y_t\right]&>\sum_{y,x} (P_\eta(y|x)-\overline P(y|x))^2\overline P(x)^2\\
     &=\sum_y\| diff_\X(y,\eta)\|^2
\end{split}
\end{equation}
To finish proving the result, let $\{\eta_i\}_{i\in\N}$ be a sequence such that 
\begin{equation}\sum_y\| diff_\X(y,\eta_i)\|^2\xrightarrow[i\to \infty]{ }0
\end{equation}
since every term is positive, then for every $y\in\Y$
\begin{equation}
\| diff_\X(y,\eta_i)\|^2\xrightarrow[i\to \infty]{ }0
\end{equation}
implying that $P_{\eta_i}(y|x)-\overline P (y|x) \xrightarrow[i\to \infty]{ }0$ for all $x,y$ and that
\begin{equation}
l(\eta_i)-l(\overline\eta) \xrightarrow[i\to \infty]{ }0
\end{equation}
Hence it's proven 
\begin{equation} (\forall\delta>0)\hspace{0.5cm} \inf_{ l(\eta)-l(\overline \eta)>\delta}\sum_y\| diff_\X(y,\eta_i)\|^2>0
\end{equation}
and therefore, after equation \ref{eq:c3_proof_key}, condition \textbf{C.3} holds.

\end{proof}
\section{Conclusion and future work}
Natural gradient based algorithms behave erratically when tested in practical problems. However, as CSNGD shows, these kind of algorithms may stabilize once convergence is guaranteed. With this in mind, we defined DSNGD, which approximates the natural gradient at each step and whose convergence in the discrete case can be proved. To that end, we stated and proved a general result showing the convergence of interior half-space gradient approximations. Furthermore, we point out that this convergence result may prove the convergence of more general algorithms, since it doesn't require the expectation of the update's direction to factor as a symmetric positive-definite matrix and the gradient.

This paper concentrates on the theoretical aspects of DSNGD. We are currently working on a flexible implementation of the algorithm that can be easily set up for different LEF linked to the conditional distributions $P_\eta (\X\mid\Y)$, including several commonly used discrete, continuous and multivariate distributions such as the normal, Poisson and exponential. Moreover, DSNGD can potentially be used in high dimensional scenarios due to its low computational complexity. The benefits of approximating the natural gradient are specially promising in this case, since the parameter space is potentially twisted and using metric information can be crucial for an algorithm's good performance. In preliminary empirical studies we are observing how it increasingly outperforms SGD as the manifold dimension grows larger. We plan to compare DSNGD against the most effective algorithms nowadays, in order to expose its weaknesses and reveal its strengths.

For the more theoretical part, in the future we plan to study the convergence of continuous and mixed DSNGD, that is when $\X$ is continuous, and in cases where $\X=(\X_d,\X_c)$ is divided into a discrete and a continuous part. 
\section{Declarations}
\subsection{Funding} 
The work has been funded by EU Horizon 2020 under grant agreements 872944 (Crowd4SDG) and 825619 (Humane-AI-net), and by the Spanish Ministry of Science and Innovation through the CI-SUSTAIN project (PID2019-104156GB-I00).
\subsection{Conflicts of interest/Competing interests}
On behalf of all authors, the corresponding author states that there is no conflict of interest.
\subsection{Availability of data and material}
Not applicable
\subsection{Code availability}
Not applicable
\subsection{Authors' contributions}
The authors have contributed equally to this work.
%
%


\begin{appendices}
\appendixpage
\addappheadtotoc
\section{Proof of Proposition \ref{propo:model}}\label{proof:lexyf_parametrization}
\begin{proof}
Prove first that if the logg-odds ratio of $P(\Y\mid\X)$ is an affine function of $\X$ then the joint distribution $P(\Y,\X)$ belongs to LEF.  

According to theorem $2$ in \cite{banerjee_analysis_2007}, assume that $P(\X\mid \Y=i)$ belongs to the same LEF for all $i\in\Y$. Also, since $\Y$ is discrete and finite, $P(\Y)$ is a categorical distribution and hence, it belongs to LEF. This means that there exist parameters $\overline{\alpha}\in\R^{s-1}$ and $\overline{\theta_i}\in\R^{t}$ for all $i\in\Y$ such that  
\begin{equation}\label{eq:pxgy}
\begin{split}
    P_{\overline{\alpha}} (\Y = i) =& \frac{\exp{S(i)^\intercal \overline{\alpha} }}{\sum_y \exp{S(y)^\intercal \overline{\alpha} }}\\
     P_{\overline{\theta_i}} (x\mid \Y = i) =&\frac{\exp{T(x)^\intercal \overline{\theta_i} }}{\int_x \exp{T(x)^\intercal \overline{\theta_i}}}\\
    \end{split}
\end{equation}
where $S$ and $T$ are sufficient statistics of $\Y$ and $\X$  respectively. If $\overline{\theta}$ is the matrix having $\overline{\theta_i}$ as $i$-th row, name $\overline{\eta}=(\overline{\alpha},\overline{\theta})$ and write
\begin{equation}
\begin{split}
    P_{\overline{\eta}}(\Y=i,x)=&P_{\overline{\alpha}}(\Y=i)P_{\overline{\theta_i}}(x\mid \Y=i)\\
    =&\frac{\exp{S(i)^\intercal \overline{\alpha} }}{\sum_y \exp{S(y)^\intercal \overline{\alpha}}}\frac{\exp{T(x)^\intercal \overline{\theta_i}}}{\int_{x} \exp{T(x)^\intercal \overline{\theta_i}}}\\
    =&\frac{\exp{S(i)^\intercal \overline{\alpha} +T(x)^\intercal \overline{\theta_i} }}{\sum_y \exp{S(y)^\intercal \overline{\alpha} }\int_{x} \exp{T(x)^\intercal \overline{\theta_i}}}
    \end{split}
\end{equation}

To prove the result, it is enough to find a change of variables from $\overline{\eta}=(\overline{\alpha},\overline{\theta})$ to $\eta=(\alpha,\beta)$ satisfying $P_{\overline{\eta}}(y,x) =P_{\eta}(y,x)$ where
\begin{equation}\label{eq:exp_fam}
    P_{\eta}(\Y=i,x)=\frac{\exp{S(i)^\intercal\alpha+T(x)^\intercal\beta_i }}{\int_x \sum_y \exp{S(y)^\intercal\alpha+T(x)^\intercal\beta_y }}
\end{equation} since $\eta$ is the natural parametrization of a LEF.

In particular, the change of variables has to satisfy that  $P_{\overline{\eta}}(x\mid \Y = i) =P_{\eta}(x\mid \Y=i)$ and $P_{\overline{\eta}}(y) =P_{\eta}(y)$. Start with the conditional probability and observe that
\begin{equation}
\begin{split}
    P_{\eta}(x\mid \Y=i)=\frac{P_{\eta}( \Y=i,x)}{\int_{x} P_{\eta}( \Y=i,x)}=&\frac{\exp{S(i)^\intercal \alpha + T(x)^\intercal \beta_i}}{\int_{x} \exp{S(i)^\intercal \alpha + T(x)^\intercal \beta_i}}\\
    =&\frac{\exp{T(x)^\intercal \beta_i}}{\int_{x}\exp{ T(x)^\intercal \beta_i}}
    \end{split}
\end{equation}

Last equation matches exactly with equation $\ref{eq:pxgy}$ by just setting $\beta = \overline{\theta}$. To complete the change of variables continue by matching $P_{\overline{\eta}}(y) =P_{\eta}(y)$.
\begin{equation}
\begin{split}
P_\eta (\Y=i) =& \frac{\int_x \exp{S(i)^\intercal \alpha +T(x)^\intercal \beta_i}}{\sum_j\int_x\exp{S(j)^\intercal \alpha +T(x)^\intercal \beta_j}}\\
=& \frac{\exp{(S(i)^\intercal \alpha)}\int_x \exp{T(x)^\intercal \beta_i}}{\sum_j\exp{(S(j)^\intercal \alpha)} \int_x\exp{T(x)^\intercal \beta_j}}\\
=& \frac{\exp{S(i)^\intercal \alpha+\log A_i}}{\sum_j\exp{S(j)^\intercal \alpha +\log A_j}}
\end{split}
\end{equation}
where $A_i = \int_x \exp{T(x)^\intercal \beta_i}$. Last equation must coincide with equation \ref{eq:pxgy}. That is
\begin{equation}\label{eq:pyinnatural}
\begin{split}
    P_\eta (\Y=i) = P_{\overline{\eta}}(\Y=i) \iff
    \frac{\exp{S(i)^\intercal \alpha+\log A_i}}{\sum_j\exp{S(j)^\intercal \alpha +\log A_j}} =
    \frac{\exp{S(i)^\intercal \overline{\alpha} }}{\sum_y \exp{S(y)^\intercal \overline{\alpha} }}
    \end{split}
\end{equation}
To simplify, assume $S$ is canonical. That is $S(i)=e_i$ is the $i$-th canonical vector for all $i\neq s$ and $S(s) = 0\in\R^{s-1}$. Note that it is enough to prove that there exists a $\mu\in\R$ such that
\begin{align}\label{eq:findmu}
    S(i)^\intercal\alpha + \log A_i -\mu= S(i)^\intercal \overline{\alpha} ,\hspace{1cm} \forall i\in\Y
\end{align}
because as a consequence, equation $\ref{eq:pyinnatural}$ clearly holds. In our case, it is $S(i)^\intercal\alpha = \alpha_i$ when $i\neq s$ and $S(i)^\intercal\alpha = 0$, and therefore the solution is
\begin{equation}
\begin{split}
   \alpha + \begin{pmatrix}
    \log A_1\\
    \vdots\\
    \log A_{s-1}
    \end{pmatrix}- 
     \begin{pmatrix}
    1\\
    \vdots\\
    1
    \end{pmatrix}
    \cdot \mu=& \overline{\alpha} \\
    \mu =& \log A_s
    \end{split}
\end{equation}
and the proof is completed when $S$ is canonical. 

Prove now the result for a general sufficient statistic $S$. Equation \ref{eq:findmu} describes the  below linear equations system
\begin{equation} 
\begin{split}\textbf{S}\alpha + \begin{pmatrix}
    \log A_1\\
    \vdots\\
    \log A_{s-1}
    \end{pmatrix}- 
     \begin{pmatrix}
    1\\
    \vdots\\
    1
    \end{pmatrix}
    \cdot \mu
    &= 
    \textbf{S}\overline{\alpha} \\
     S(s)^\intercal\alpha + \log A_s -\mu &= S(s)^\intercal \overline{\alpha}
     \end{split}
\end{equation}
where $\textbf{S}$ is the matrix having $S(1),...,S(s-1)$ as rows.  Since $S$ is a sufficient statistic, assume without loss of generality that $S(1),...,S(s-1)$ are linearly independent vectors, and then $\textbf{S}$ is invertible. Finally, it is easy to check that the change of variables is
\begin{equation}
\begin{split}
 \alpha + \textbf{S}^{-1}\left(\begin{pmatrix}
    \log A_1\\
    \vdots\\
    \log A_{s-1}
    \end{pmatrix}-
     \begin{pmatrix}
    1\\
    \vdots\\
    1
    \end{pmatrix}
    \cdot \mu\right)
    &= 
    \overline{\alpha}\\
    \mu&=\frac{S(s)^\intercal \textbf{S}^{-1}\begin{pmatrix}
    \log A_1\\
    \vdots\\
    \log A_{s-1}
    \end{pmatrix}
    - \log A_s}{S(s)^\intercal \textbf{S}^{-1}\begin{pmatrix}
    1\\
    \vdots\\
    1
    \end{pmatrix}
    - 1}
    \end{split}
\end{equation}

The converse implication is straightforward. Assuming that $P(\X,\Y)$ belongs to LEF, and therefore assuming equation \ref{eq:exp_fam}, start by expressing the conditional probability distributions of $\Y$ given $\X$ in $\eta$.
\begin{equation}
\begin{split}
    P_{\eta}(y\mid x) = &\frac{P_{\eta}(y,x)}{\sum_y P_{\eta}(y,x)}\\=&
    \frac{\exp S(y)^\intercal\alpha+T(x)^\intercal\beta_y}{\sum_y \exp S(y)^\intercal\alpha+T(x)^\intercal\beta_y}
    \end{split}
\end{equation}
and compute the log-odds ratio
\begin{equation}
    \begin{split}
    \log{\frac{P_\eta(x\mid \Y=k)}{P_\eta(x\mid \Y=h)}} =& S(k)^\intercal\alpha+T(x)^\intercal\beta_k -(S(h)^\intercal\alpha+T(x)^\intercal\beta_h)\\
    =& (S(k)-S(h))^\intercal \alpha +  T(x)^\intercal(\beta_k-\beta_h)
    \end{split}
\end{equation}
which is clearly an affine function of features $\X$.
\end{proof}

\section{Proof of Proposition \ref{propo:nat_grad_dsngd}}\label{proof:nat_grad_dsngd}
\begin{proof}

First, claim that

\begin{equation}\label{eq:log-loss_grad}
    {\nabla}l(\eta,x,y) = {\nabla} h(x,\eta)\cdot (q_\Y(x,P_{\eta})-e_s(y))
\end{equation}
Indeed,
\begin{equation}
    \begin{split}
     {\nabla}  \log P_{\eta} (y\mid x) =& {\nabla} \log P_{\eta}(y, x)-{\nabla}\log \sum_y P_{\eta}(y, x)\\
    =& {\nabla} \log P_{\eta}(y, x)-\frac{\sum_y  {\nabla} P_{\eta}(y, x)}{\sum_y P_{\eta}(y,x)}\\
    =& {\nabla} \log P_{\eta}(y, x)-\frac{\sum_y P_{\eta}(y, x) {\nabla} \log P_{\eta}(y, x)}{\sum_y P_{\eta}(y,x)}\\
    =&{\nabla} \log P_{\eta}(y, x)-\sum_y P_{\eta}(y\mid x) {\nabla} \log P_{\eta}(y, x)\\
    =&{\nabla} h(y, x,\eta)-\E_{\Y\mid x}[{\nabla} h(y,x,\eta)]\label{eq:grad_expect}   
    \end{split}
\end{equation}

where $h(y,x,\eta) =\log P_{\eta}(y, x)$.
Observe we can rewrite equation \ref{eq:grad_expect} as;
    \begin{equation}
    {\nabla} \log P_{\eta} (y\mid x) = -{\nabla} h( x,\eta)\cdot (q(x,\eta)-e_s(i))
\end{equation}
where $ h(x,\eta)=(h(1,x,\eta),...,h(s,x,\eta))$ implying the claim. From equation~\ref{eq:log-loss_grad} observe that

\begin{equation}
    \widetilde{\nabla}l(\eta,x,y) = \widetilde{\nabla} h(x,\eta)\cdot (q_\Y(x,P_{\eta})-e_s(y))
    \end{equation}
Finally, since the log-loss is defined in a DFM, then use previous equation and equation~\ref{eq:natural_gradient_is_gradient} to finish the proof
\end{proof}

\section{Proof of Proposition \ref{theorem:h}}\label{proof:h}
\begin{proof}
To simplify, break $\nabla_{\eta^*} = (\nabla_{\alpha^*},\nabla_{\beta_1^*},...,\nabla_{\beta_s^*})$ and then it's clear that
\begin{align}\nabla h(x,\eta^*)=
    \begin{pmatrix}
     \nabla_{\alpha^*} h(x,\eta^*)\\
     \nabla_{\beta_1^*} h(x,\eta^*)  \\
     \vdots\\
     \nabla_{\beta_s^*} h(x,\eta^*)  
    \end{pmatrix}
\end{align}
Start with $\nabla_{\alpha^*} h(x,\eta^*)$ expression. Observe that $i$-th column of $\nabla_{\alpha^*} h(x,\eta^*)$ is
\begin{equation}
\begin{split}
\nabla_{\alpha^*} \log P_{\eta^*}(\Y= i,x) =& \nabla_{\alpha^*}\log P_{\alpha^*}(\Y=i)+\nabla_{\alpha^*}\log P_{\alpha^*}(x\mid \Y=i)\\
=& \nabla_{\alpha^*}\log P_{\alpha^*}(\Y=i)+d_{\alpha^*} \theta_i \nabla_{\theta_i} \log P_{\theta_i}(x\mid \Y=i)
\end{split}
\end{equation}
where in last step the chain rule is applied and $d_{\alpha^*} \theta_i$ stands for the Jacobian of $\theta_i$ with respect to $\alpha^*$. 

Assume the canonical parametrization is used, then according to equation \ref{eq:expect_parametriztion} write
\begin{equation}\label{eq:y_probabilities}
    \alpha^* =\begin{pmatrix}
    P_{\eta^*}(\Y=1)\\
    \vdots\\
    P_{\eta^*}(\Y=s-1)
    \end{pmatrix}
\end{equation}

From equations $\ref{eq:y_probabilities}$ and $\ref{eq:expect_param_x|y}$ obtain
 \begin{equation}
 \begin{split}
      \nabla_{\alpha^*}\log P_{\alpha^*}(\Y=i)=&\frac{1}{P_{\alpha^*}(\Y=i)}
      \left\{\begin{array}{cc}
           e_{s-1}(i) & i\neq s \\
          (-\textbf{1}) & i=s
      \end{array}\right.\\
       d_{\alpha^*}\theta_i=&\frac{-1}{P_{\alpha^*}(\Y=i)}
      \left\{\begin{array}{cc}
          e_{s-1}(i) \cdot \theta_i^\intercal & i\neq s \\
          (-\textbf{1}) \cdot \theta_i^\intercal& i=s
      \end{array}\right.
      \end{split}
 \end{equation}
 where $e_{s-1}(i)$ is the $i$-th canonical $s-1$ dimensional vector and $\textbf{1}=
 \begin{pmatrix}
 1\\
 \vdots\\
 1
 \end{pmatrix}
 $. From here deduce,
 
\begin{equation}
\begin{split}
    \nabla_{\alpha^*}h(x,\eta^*)=&
    K_{(s-1)\times s}
    \cdot diag(
    d(x,1,\zeta^*),...,d(x,s,\zeta^*))
    \\
    d(x,y,\zeta^*)=&\frac{1-\theta_y^\intercal \nabla_{\theta_y}\log P(x\mid y)}{P_{\zeta^*}(y)}
    \end{split}
\end{equation}
The part  $\nabla_{\beta_k^*} h(x,\eta^*)$ follows the same steps. Observe that $i$-th column of $\nabla_{\beta_y^*} h(x,\eta^*)$ is
\begin{equation}
\begin{split}
\nabla_{\beta_y^*} \log P_{\eta^*}(\Y= i,x) =& \nabla_{\beta_y^*}\log P_{\beta_y^*}(x\mid \Y)\\
=& d_{\beta_y^*} \theta_i \nabla_{\theta_i} \log P_{\theta_i}(x\mid \Y=i)\\
=& \begin{cases}
0 & y\neq i\\
\frac{ \nabla_{\theta_i} \log P_{\theta_i}(x\mid \Y=i)}{P_{\zeta^*}(y)} & y=i
\end{cases}
\end{split}
\end{equation}
and therefore the claim is proved 
\end{proof}

\section{Proof of Proposition \ref{theorem:h_discrete}}\label{proof:h_discrete}

\begin{proof}
Compute $\nabla_{\theta_y} \log P_{\theta_y}(x|y) $ and proposition~\ref{theorem:h} finishes the proof. Parameters $\theta_y= (\theta_{y,1},...,\theta_{y,m-1})$ are the expectation parameters of the probability distribution $ P_{\theta_y}(x|y)$ which belongs to LEF.

Recall that the canonical parametrization is taken and by equation~\ref{eq:expect_param_x|y} deduce
\begin{equation}
    \begin{split}
        P_{\theta_y}(x|y) = 
        \begin{cases}
        \theta_{y,x} & x\neq m\\
        1-\sum_j \theta_{y,j} & x=m
        \end{cases}
    \end{split}
\end{equation}

which clearly implies
\begin{equation}
    \begin{split}
        \nabla_{\theta_y} \log P_{\theta_y}(x|y)= 
         \frac{1}{ P_{\theta_y}(x|y)}\begin{cases}
        e_{m-1}(x)& x\neq m\\
        \textbf{-1}_{m-1}  & x=m
        \end{cases}
    \end{split}
\end{equation}

Finally observe
\begin{equation}
    d(x,y,\zeta^*)=\frac{1-\theta_y^\intercal \nabla_{\theta_y}\log P_{\theta_y}(x\mid y)}{P_{\zeta^*}(y)}=
    \begin{cases}
    0 & x\neq m\\
    \frac{1}{P_{\theta_y}(x|y)P_{\zeta^*}(y)} & x=m
    \end{cases}
\end{equation}

Substitute the computations in proposition~\ref{theorem:h} to finish the proof.
\end{proof}

\section{Proof of Proposition \ref{propo:nat_grad_complexity}}\label{proof:natural_gradient_complexity}

\begin{proof}
Let $A$ be the cost of computing $\nabla_{\theta_k}\log P_{\theta_k}(x\mid y)$. Prove first the next claim: the number of operations required to compute  $\nabla_{\zeta^*}h(x,\zeta^*)$ is 
\begin{equation}
    s\cdot (A  + 3 t+2)-1
    \end{equation}and hence $O(s\cdot(A+t))$.

Indeed, terms $\nabla_{\theta_k}\log P_{\theta_k}(x\mid y)$  for every $y\in\Y$ need $s\cdot A$ operations. The cost of computing $P_{\zeta^*}(y)$ for every $y\in\Y$ is $s-1$ according to equation~\ref{eq:expect_parametriztion} (only the term $P_{\zeta^*}(s)$ requires operations). Obtain term $d(x,y,\zeta^*)$ after $2t+1$ operations ( $2t-1$ for the scalar product of vectors, $1$ for the subtraction in the numerator and $1$ last operation for the division). Since this needs to be done for every $y\in\Y$ then $d(x,1,\zeta^*),...,d(x,s,\zeta^*)$ is known with $s\cdot (2t+1)$ operations. Now, $\nabla_{\alpha^*} h(x,\zeta^*)$ is obtained with the product of matrices $M$ (which is almost the identity matrix) and a diagonal matrix, which does not require any operation (it is just a transformation). Finally $\nabla_{\beta^*_k} h(x,\zeta^*)$ demands for $t$ divisions for every $y\in\Y$, and therefore for $s\cdot t$ operations. Hence, the claim is proved.

To previous analysis, add the costs represented by equation~\ref{eq:nat_grad_dsngd}. That is, analyze the costs of computing $q_\Y(x,P_\eta)$ and then the products shown in that equation.

The vector $q_\Y(x,P_\eta)$ consist on computing $P_\eta(y\mid x)$ for every $y\in\Y$. Using equation~\ref{eq:conditionals}, $q_\Y(x,P_\eta)$ needs $2t+1$ for scalar products $T(x)^\intercal \beta_y$, $1$ subtraction in $S(y)^\intercal\alpha -T(x)^\intercal \beta_y$ (recall that $S$  statistic is canonical), then $1$ exponentiation and finally $1$ division.  This is done for every $y\in \Y$. The denominator is the same for every $y$ so it can be computed just once with $s-1$ sums. The total is
\begin{equation}
    2ts+5s-1
\end{equation} operations.

Finally, the  operations described in equation~\ref{eq:nat_grad_dsngd} are $1$ for $(q_\Y(x,P_{\eta})-e_s(y))$, $s-1$ for $\nabla_{\alpha^*} h(x,\zeta^*)\cdot \nabla(q_\Y(x,P_{\eta})-e_s(y))$ product and $t$ operations for $\nabla_{\beta^*_y} h(x,\zeta^*)\cdot \nabla(q_\Y(x,P_{\eta})-e_s(y))$, this last one needs to be done for every $y\in\Y$. The total operations for this block it is then 
\begin{equation}
    s+s\cdot t
\end{equation}

To conclude the proof, the total operations needed is 
\begin{equation}
    s\cdot(A+6t+8) -2
\end{equation}
and the complexity order is $O(s\cdot(A+t))$

\end{proof}

\section{Proof of Theorem \ref{theorem:convergence}}\label{proof:convergence_generalization}

\begin{proof}
The proof uses Robbins-Siegmund theorem as key tool. Steps taken are closely inspired by those taken in the proof of Theorem 3.2 in \cite{sunehag_variable_2009}. 

Compute Taylor' second order approximation of $l(\eta_{t+1})$, and after condition \textbf{C.2} apply Taylor's inequality 
\begin{equation}
    l(\eta_{t+1}) \leq  l(\eta_t)-\gamma_t \nabla l(\eta_t)^T  Y_t+\gamma_t^2 K \| Y_t\|^2
\end{equation}
Therefore, applying the expectation conditioned to information at time $t$ obtain
\begin{equation}\label{eq:exp_taylor}
\E_t[ l(\eta_{t+1})] \leq
     l(\eta_t)-\gamma_t \nabla l(\eta_t)^T \E_t [ Y_t]+\gamma_t^2 K \E_t \| Y_t\|^2
\end{equation}
 Use bound of \textbf{C.4} to third term of right hand side
\begin{equation}
\E_t [l(\eta_{t+1})] \leq l(\eta_t) -\gamma_t \nabla l(\eta_t)^T \E_t [ Y_t] + \gamma_t^2 K (A+Bl(\eta_t))
\end{equation}
Finally, substitute $U_t = l(\eta_t)$ and arrange terms to match with equation~\ref{eq:Robbins-Siegmund}
\begin{equation}\label{eq:R-S_match}
    \E_t [U_{t+1}]\leq (1+B\gamma_t^2 K)U_t -\gamma_t \nabla l(\eta_t)^T \E_t [ Y_t]+ \gamma_t^2 K A
\end{equation}
Note that theorem \ref{theorem:Robbins-Siegmund} conditions are satisfied, since condition \textbf{C.6} implies $\sum_t\beta_t =\sum_t BK\gamma^2=BK\sum_t \gamma^2<\infty$ and $\sum_t\epsilon_t =\sum_t KA\gamma^2<\infty$. Hence, Robbins-Siegmund theorem ensures that $U_{t} = l(\eta_t)$ converges almost surely to a random variable and 
\begin{equation}\label{eq:finite_sum}
\sum_t \zeta_t=\sum_t \gamma_t \nabla l(\eta_t)^T \E_t [ Y_t]<\infty
\end{equation} 
Now prove that $\lim_t l(\eta_t) = l(\overline\eta)$. If $l(\eta_t)$ converges to some different random variable, condition \textbf{C.3}, second condition of \textbf{C.6} and equation \ref{eq:finite_sum} lead to a contradiction. Indeed, if $\lim_t l(\eta_t)=v\neq l(\overline\eta)$, use condition \textbf{C.3} and deduce that for a fixed $0<\delta<v-l(\overline\eta)$ there exists an $N$ large enough and $\epsilon>0$ such that
\begin{equation}
    \nabla l(\eta_t)^T\E_t \left[ Y_t\right]\geq \epsilon
\end{equation}
for all $t>N$. Therefore, equation \ref{eq:finite_sum} becomes
\begin{equation}
\begin{split}
    \sum_t \gamma_t \nabla l(\eta_t)^T \E_t [ Y_t] &= \sum_t ^N \gamma_t \nabla l(\eta_t)^T \E_t [ Y_t] + \sum_{t>N}\gamma_t \nabla l(\eta_t)^T \E_t [ Y_t]\\
     &\geq\sum_t ^N \gamma_t \nabla l(\eta_t)^T \E_t [ Y_t] + \sum_{t>N} \epsilon\gamma_t\\
     &\geq\epsilon\sum_{t>N} \gamma_t
    \end{split}
\end{equation}
Second condition in \textbf{C.6} applied to right hand side of above equation assures that
\begin{equation}
    \sum_t \gamma_t \nabla l(\eta_t)^T \E_t [ Y_t] = \infty
\end{equation}
which contradicts equation \ref{eq:finite_sum}.

Finally, it is only possible that $\lim_t l(\eta_t) = l(\overline\eta)$ almost surely as we wanted to prove.
\end{proof}
\section{Proof of condition \textbf{C.2} in Theorem \ref{theorem:dsngd_convergence}}\label{proof:C2}
\begin{proof}
Compute the hessian of 

\begin{equation}
    l(\eta)=\sum_{x,y} l(\eta, y, x) \overline P(y,x)
\end{equation}

The gradient of $ l(\eta, y, x)$ is
\begin{equation}\label{eq:grad}
    \begin{split}
        \nabla_\alpha l(\eta, y,x) &=S\cdot (q_\Y(x)-e_s(y))\\
        \nabla_{\beta_{y'}} l(\eta, y, x) &=(q_\Y( x)_{y'}-\delta_{y=y'})\cdot T(x)
    \end{split}
\end{equation}
where $S$ is the matrix having $S(i)$ as $i$-th column for $i\in\Y$. Therefore, the hessian is
\begin{equation}\label{eq:l_hessian}
    \begin{split}
        \nabla^2_\alpha l(\eta, y,x) &=S\cdot( diag(q_\Y(x)) -q_\Y(x)\cdot q_\Y(x)^\intercal)\cdot S^\intercal\\
        \nabla_{\beta_{y_2}}\nabla_{\beta_{y_1}} l(\eta,y,x) &=\nabla_{\beta_{y_2}}q_\Y(x)_{y_1}\cdot T(x)\\
        &= -T(x)\cdot T(x)^\intercal q_\Y(x)_{y_1}(q_\Y(x)_{y_2}-\delta_{y_1=y_2})\\
        \nabla_{\alpha}\nabla_{\beta_{y'}} l(\eta,y,x) &=\nabla_{\alpha}q_\Y(x)_{y'}\cdot T(x)\\
        &=T(x)\cdot (q_\Y(x)-e_s(y')^\intercal \cdot S^\intercal
    \end{split}
\end{equation}

Observe how all matrices in equation \ref{eq:l_hessian} have their elements bounded once $S$ and $T$ statistics are fixed, since $\|q_\Y(x)\|\leq 1$. Therefore
\begin{equation}
     \|\nabla^2l(\eta,y,x)\|\leq 2 K_{x,y}
\end{equation}
for some positive numbers $K_{x,y}$. Define $K=\max_{x,y} K_{x,y}$. then  finally
\begin{equation}
    \begin{split}
        \|\nabla^2_\eta l(\eta)\|=&\|\nabla^2\sum_{y,x}l(\eta,y,x)\cdot \overline P(x,y)\|\\
        =&\|\sum_{y,x} \nabla^2l(\eta,y,x)\cdot \overline P(x,y)\|\\
        \leq &\sum_{y,x} \|\nabla^2l(\eta,y,x)\|\cdot \overline P(x,y)\\
        \leq&\sum_{y,x} 2 \cdot K_{x,y}\cdot \overline P(x,y)\\
        \leq&2\cdot K\sum_{y,x} \overline P(x,y)\\
        =&2\cdot K
    \end{split}
\end{equation}
\end{proof}

\section{Proof of condition \textbf{C.4} in Theorem \ref{theorem:dsngd_convergence}}\label{proof:C4}

\begin{proof}

Observe that for any $\epsilon$ and $t$ large enough there exists $A_{x_t}$ such that
\begin{equation}
\begin{split}
    \| Y_t\|^2&= {(q_\Y(x_t)-e_s({y_t}))^\intercal\cdot h(x_t,\zeta_t^*)^\intercal h(x_t,\zeta_t^*)\cdot(q_\Y(x_t)-e_s({y_t}))}\\
    &\leq A_{x_t}\|q_\Y(x_t)-e_s({y_t})\|^2
    \end{split}
\end{equation}
where 
\begin{equation}
A_{x_t}\geq\|h(x_t,\zeta_t^*)^\intercal h(x_t,\zeta_t^*)\|+\epsilon
\end{equation}
This is because $\zeta_t$ converges and because of theorem~\ref{theorem:h_discrete}. Now
\begin{equation}
\begin{split}
    \|q_\Y(x_t)-e_s({y_t})\|^2&=1-2P_{\eta_t}(y_t|x_t)+\sum_y{P_{\eta_t}(y|x_t)^2}\\
    &\leq s+1 
    \end{split}
\end{equation}
therefore $\| Y_t\|^2\leq A_{x_t}(s+1)$ and
\begin{equation}
\begin{split}
    \E\| Y_t\|^2 &\leq \E[ A_{x_t}(s+1)]\\
    &\leq A'(s+1)=A
\end{split}
\end{equation}
where $A'=\max_x A_{x}$ and then condition \textbf{C.4} holds.
\end{proof}
\end{appendices}

%
%

\bibliography{main}   

\begin{thebibliography}{10}
\providecommand{\url}[1]{{#1}}
\providecommand{\urlprefix}{URL }
\expandafter\ifx\csname urlstyle\endcsname\relax
  \providecommand{\doi}[1]{DOI~\discretionary{}{}{}#1}\else
  \providecommand{\doi}{DOI~\discretionary{}{}{}\begingroup
  \urlstyle{rm}\Url}\fi

\bibitem{nemirovski_robust_2009}
Nemirovski, A., Juditsky, A., Lan, G., Shapiro, A.: Robust {Stochastic}
  {Approximation} {Approach} to {Stochastic} {Programming}.
\newblock SIAM Journal on Optimization \textbf{19}(4), 1574--1609 (2009).
\newblock \doi{10.1137/070704277}.
\newblock \urlprefix\url{https://epubs.siam.org/doi/abs/10.1137/070704277}.
\newblock Publisher: Society for Industrial and Applied Mathematics

\bibitem{hu_brief_2020}
Hu, J., Liu, X., Wen, Z.W., Yuan, Y.X.: A {Brief} {Introduction} to {Manifold}
  {Optimization}.
\newblock Journal of the Operations Research Society of China \textbf{8}(2),
  199--248 (2020).
\newblock \doi{10.1007/s40305-020-00295-9}.
\newblock \urlprefix\url{https://doi.org/10.1007/s40305-020-00295-9}

\bibitem{amari_natural_1998}
Amari, S.i.: Natural {Gradient} {Works} {Efficiently} in {Learning}.
\newblock Neural Computation \textbf{276}, 251--276 (1998)

\bibitem{borja_jesus_2019}
S\'anchez-L\'opez, B., Cerquides, J.: Convergent stochastic almost natural
  gradient descent.
\newblock Artificial Intelligence Research and Development- Proceedings of the
  22nd International Conference of the Catalan Association for Artificial
  Intelligence \textbf{319}, 54--63 (2019)

\bibitem{sunehag_variable_2009}
Sunehag, P., Trumpf, J., Vishwanathan, S.V.N., Schraudolph, N.: Variable
  {Metric} {Stochastic} {Approximation} {Theory}.
\newblock In: Artificial {Intelligence} and {Statistics}, pp. 560--566 (2009).
\newblock \urlprefix\url{http://proceedings.mlr.press/v5/sunehag09a.html}

\bibitem{li_spectral_2012}
{Li}, J., {Bioucas-Dias}, J.M., {Plaza}, A.: Spectral–spatial hyperspectral
  image segmentation using subspace multinomial logistic regression and markov
  random fields.
\newblock IEEE Transactions on Geoscience and Remote Sensing \textbf{50}(3),
  809--823 (2012).
\newblock \doi{10.1109/TGRS.2011.2162649}

\bibitem{covington_deep_nn_2016}
Covington, P., Adams, J., Sargin, E.: Deep neural networks for youtube
  recommendations.
\newblock In: Proceedings of the 10th ACM Conference on Recommender Systems,
  RecSys '16, p. 191–198. Association for Computing Machinery, New York, NY,
  USA (2016).
\newblock \doi{10.1145/2959100.2959190}.
\newblock \urlprefix\url{https://doi.org/10.1145/2959100.2959190}

\bibitem{daniels1997hierarchical}
Daniels, M.J., Gatsonis, C.: Hierarchical polytomous regression models with
  applications to health services research.
\newblock Statistics in Medicine \textbf{16}(20), 2311--2325 (1997)

\bibitem{bull2007confidence}
Bull, S.B., Lewinger, J.P., Lee, S.S.: Confidence intervals for multinomial
  logistic regression in sparse data.
\newblock Statistics in Medicine \textbf{26}(4), 903--918 (2007)

\bibitem{biesheuvel2008polytomous}
Biesheuvel, C., Vergouwe, Y., Steyerberg, E., Grobbee, D., Moons, K.:
  Polytomous logistic regression analysis could be applied more often in
  diagnostic research.
\newblock Journal of clinical epidemiology \textbf{61}(2), 125--134 (2008)

\bibitem{leppink2020multicategory}
Leppink, J.: Multicategory nominal choices.
\newblock In: The Art of Modelling the Learning Process, pp. 103--110. Springer
  (2020)

\bibitem{ben-akiva_discrete_1985}
Ben-Akiva, M.E., Lerman, S.R., Lerman, S.R.: Discrete {Choice} {Analysis}:
  {Theory} and {Application} to {Travel} {Demand}.
\newblock MIT Press (1985).
\newblock Google-Books-ID: oLC6ZYPs9UoC

\bibitem{tadei_recent_2018}
Tadei, R., Perboli, G., Manerba, D.: A {Recent} {Approach} to {Derive} the
  {Multinomial} {Logit} {Model} for {Choice} {Probability}.
\newblock In: P.~Daniele, L.~Scrimali (eds.) New {Trends} in {Emerging}
  {Complex} {Real} {Life} {Problems}: {ODS}, {Taormina}, {Italy}, {September}
  10–13, 2018, {AIRO} {Springer} {Series}, pp. 473--481. Springer
  International Publishing, Cham (2018)

\bibitem{NIPS2008_077e29b1}
Nock, R., Nielsen, F.: On the efficient minimization of classification
  calibrated surrogates.
\newblock In: D.~Koller, D.~Schuurmans, Y.~Bengio, L.~Bottou (eds.) Advances in
  Neural Information Processing Systems, vol.~21, pp. 1201--1208. Curran
  Associates, Inc. (2009)

\bibitem{reidCompositeBinaryLosses2010}
Reid, M.D., Williamson, R.C.: Composite {{Binary Losses}}.
\newblock Journal of Machine Learning Research \textbf{11}(83), 2387--2422
  (2010)

\bibitem{vernetCompositeMulticlassLosses2011}
Vernet, E., Reid, M.D., Williamson, R.C.: Composite {{Multiclass Losses}}.
\newblock In: J.~{Shawe-Taylor}, R.S. Zemel, P.L. Bartlett, F.~Pereira, K.Q.
  Weinberger (eds.) Advances in {{Neural Information Processing Systems}} 24,
  pp. 1224--1232. {Curran Associates, Inc.} (2011)

\bibitem{nockSupervisedLearningNo2020}
Nock, R., Menon, A.K.: Supervised {{Learning}}: {{No Loss No Cry}}.
\newblock arXiv:2002.03555 [cs, stat]  (2020)

\bibitem{vapnik_principles_1991}
Vapnik, V.: Principles of risk minimization for learning theory.
\newblock In: Proceedings of the 4th {International} {Conference} on {Neural}
  {Information} {Processing} {Systems}, {NIPS}'91, pp. 831--838. Morgan
  Kaufmann Publishers Inc., San Francisco, CA, USA (1991)

\bibitem{bottou-98x}
Bottou, L.: Online algorithms and stochastic approximations.
\newblock In: D.~Saad (ed.) Online Learning and Neural Networks. Cambridge
  University Press, Cambridge, UK (1998).
\newblock \urlprefix\url{http://leon.bottou.org/papers/bottou-98x}.
\newblock Revised, oct 2012

\bibitem{dennis_numerical_1996}
Dennis, J.E., Schnabel, R.B.: Numerical methods for unconstrained optimization
  and nonlinear equations.
\newblock No.~16 in Classics in applied mathematics. SIAM, Philadelphia, Pa
  (1996).
\newblock OCLC: 845110213

\bibitem{becker_improving_1989}
Becker, S., Lecun, Y.: Improving the convergence of back-propagation learning
  with second-order methods.
\newblock In: D.~Touretzky, G.~Hinton, T.~Sejnowski (eds.) Proceedings of the
  1988 {Connectionist} {Models} {Summer} {School}, {San} {Mateo}, pp. 29--37.
  Morgan Kaufmann (1989)

\bibitem{duchi_adaptive_2011}
Duchi, J., Hazan, E., Singer, Y.: Adaptive subgradient methods for online
  learning and stochastic optimization.
\newblock J. Mach. Learn. Res. \textbf{12}(null), 2121–2159 (2011)

\bibitem{zeiler_adadelta:_2012}
Zeiler, M.D.: {ADADELTA}: {An} {Adaptive} {Learning} {Rate} {Method}.
\newblock arXiv:1212.5701 [cs]  (2012).
\newblock \urlprefix\url{http://arxiv.org/abs/1212.5701}.
\newblock ArXiv: 1212.5701

\bibitem{tieleman_lecture_2012}
Tieleman, T., Hinton, G.: Lecture 6.5-rmsprop: Divide the gradient by a running
  average of its recent magnitude.
\newblock COURSERA: Neural networks for machine learning \textbf{4}(2), 26--31
  (2012)

\bibitem{kingma_adam:_2015}
Kingma, D.P., Ba, L.J.: Adam: {A} {Method} for {Stochastic} {Optimization}.
\newblock In: Proceedings of the 3rd {International} {Conference} on {Learning}
  {Representations} ({ICLR}) (2015).
\newblock
  \urlprefix\url{https://dare.uva.nl/search?identifier=a20791d3-1aff-464a-8544-268383c33a75}

\bibitem{carmo_riemannian_2013}
Carmo, M.P.d.: Riemannian geometry, corrected at 14th printing\$h2013 edn.
\newblock Mathematics: theory \& applications. Birkhäuser, Boston Basel Berlin
  (2013)

\bibitem{murray1993differential}
Murray, M.K., Rice, J.W.: Differential geometry and statistics, vol.~48.
\newblock CRC Press (1993)

\bibitem{thomas_2014}
Thomas, P.S.: Genga: A generalization of natural gradient ascent with positive
  and negative convergence results.
\newblock 31st International Conference on Machine Learning, ICML 2014
  \textbf{5}, 3533--3541 (2014)

\bibitem{bonnabel_stochastic_2013}
Bonnabel, S.: Stochastic gradient descent on {Riemannian} manifolds.
\newblock IEEE Transactions on Automatic Control \textbf{58}(9), 2217--2229
  (2013).
\newblock \doi{10.1109/TAC.2013.2254619}.
\newblock \urlprefix\url{http://arxiv.org/abs/1111.5280}.
\newblock ArXiv: 1111.5280

\bibitem{amari_information_2016}
Amari, S.i.: Information geometry and its applications, vol. 5416.
\newblock Springer (2016)

\bibitem{nielsen_elementary_2018}
Nielsen, F.: An elementary introduction to information geometry.
\newblock arXiv:1808.08271 [cs, math, stat]  (2018).
\newblock \urlprefix\url{http://arxiv.org/abs/1808.08271}.
\newblock ArXiv: 1808.08271

\bibitem{nemirovski1983}
Nemirovski{\u\i} A~\and~Yudin, D.: Problem Complexity and Method Efficiency in
  Optimization.
\newblock A Wiley-Interscience publication. Wiley (1983).
\newblock \urlprefix\url{https://books.google.es/books?id=6ULvAAAAMAAJ}

\bibitem{masegosa_stochastic_2014}
Masegosa, A.R.: Stochastic {Discriminative} {EM}.
\newblock In: Proceedings of the {Thirtieth} {Conference} on {Uncertainty} in
  {Artificial} {Intelligence}, {UAI}'14, pp. 573--582. AUAI Press, Arlington,
  Virginia, United States (2014).
\newblock \urlprefix\url{http://dl.acm.org/citation.cfm?id=3020751.3020811}.
\newblock Event-place: Quebec City, Quebec, Canada

\bibitem{raskutti_information_2015}
Raskutti, G., Mukherjee, S.: The {Information} {Geometry} of {Mirror}
  {Descent}.
\newblock IEEE Transactions on Information Theory \textbf{61}(3), 1451--1457
  (2015).
\newblock \doi{10.1109/TIT.2015.2388583}

\bibitem{BECK2003167}
Beck, A., Teboulle, M.: Mirror descent and nonlinear projected subgradient
  methods for convex optimization.
\newblock Operations Research Letters \textbf{31}(3), 167 -- 175 (2003).
\newblock \doi{https://doi.org/10.1016/S0167-6377(02)00231-6}.
\newblock
  \urlprefix\url{http://www.sciencedirect.com/science/article/pii/S0167637702002316}

\bibitem{banerjee_analysis_2007}
Banerjee, A.: An {Analysis} of {Logistic} {Models}: {Exponential} {Family}
  {Connections} and {Online} {Performance}.
\newblock In: Proceedings of the 2007 {SIAM} {International} {Conference} on
  {Data} {Mining}, Proceedings, pp. 204--215. Society for Industrial and
  Applied Mathematics (2007).
\newblock \doi{10.1137/1.9781611972771.19}.
\newblock
  \urlprefix\url{https://epubs.siam.org/doi/abs/10.1137/1.9781611972771.19}

\bibitem{wani_1968}
Wani, J.K.: On the linear exponential family.
\newblock Mathematical Proceedings of the Cambridge Philosophical Society
  \textbf{64}(2), 481–483 (1968).
\newblock \doi{10.1017/S0305004100043097}

\bibitem{robbins_siegmund_1971}
Robbins, H., Siegmund, D.: A convergence theorem for non negative almost
  supermartingales and some applications.
\newblock In: J.S. Rustagi (ed.) Optimizing Methods in Statistics, pp. 233 --
  257. Academic Press (1971)

\end{thebibliography}
\bibliographystyle{spmpsci}      

%
%

\end{document}